\documentclass[11pt]{article}
\usepackage{amsfonts}
\usepackage{graphics}
\usepackage{indentfirst}
\usepackage{cite}
\usepackage{latexsym}
\usepackage[colorlinks, linkcolor=red]{hyperref}
\usepackage{color, xcolor}
\usepackage[paper=a4paper, left=1.6cm, right=1.6cm, top=1.8cm, bottom=1.6cm, headheight=5.5pt, footskip=0.8cm, footnotesep=0.8cm, centering, includefoot]{geometry}
\hypersetup{urlcolor=red, citecolor=red}
\usepackage{amsmath}
\allowdisplaybreaks
\usepackage{amssymb}
\usepackage[dvips]{epsfig}
\usepackage{amscd}

\newtheorem{theorem}{Theorem}[section]
\newtheorem{remark}{Remark}[section]

\newtheorem{lemma}{Lemma}[section]

\DeclareMathOperator{\divv}{div}
\DeclareMathOperator{\curl}{curl}

\makeatletter
\@addtoreset{equation}{section}
\makeatother
\makeatletter
\@addtoreset{equation}{section}
\makeatother

\title{Global well-posedness to the 2D Cauchy problem of nonhomogeneous heat conducting Navier-Stokes and magnetohydrodynamic equations with vacuum at infinity
\thanks{This research was partially supported by National Natural Science Foundation of China (Nos. 11901474, 12071359), Exceptional Young Talents Project of Chongqing Talent (No. cstc2021ycjh-bgzxm0153), and the Innovation Support Program for Chongqing Overseas Returnees (No. cx2020082).}
}

\date{}
\author{ Xin Zhong\thanks{School of Mathematics and Statistics, Southwest University, Chongqing 400715, People's Republic of China
({\tt xzhong1014@amss.ac.cn}).
}
}

\begin{document}
\maketitle

\begin{abstract}
We revisit the 2D Cauchy problem of nonhomogeneous heat conducting magnetohydrodynamic (MHD) equations in $\mathbb{R}^2$. For the initial density allowing vacuum at infinity, we derive the global existence and uniqueness of strong solutions provided that the initial density and the initial magnetic decay not too slowly at infinity. In particular, the initial data can be arbitrarily large. This improves our previous work \cite{zx2} where the initial density has non-vacuum states at infinity.
The result could also be viewed as an extension of the study in L{\"u}-Xu-Zhong \cite{lvb2} for the inhomogeneous case to the full inhomogeneous situation. The method is based on delicate spatial weighted estimates and the structural characteristic of the system under consideration. As a byproduct, we get the global existence of strong solutions to the 2D Cauchy problem for nonhomogeneous heat conducting Navier-Stokes equations with vacuum at infinity.
\end{abstract}

\textit{Key words and phrases}. Nonhomogeneous heat conducting MHD equations;  global well-posedness; 2D Cauchy problem; vacuum at infinity.

2020 \textit{Mathematics Subject Classification}. 76D05; 76W05; 76D03.

\section{Introduction and main results}

Magnetohydrodynamics is the study of the interaction of electromagnetic
fields and conducting fluids. The modeling consists of a coupling between the
Navier-Stokes equations of continuum fluid mechanics and the Maxwell equations of electromagnetism. In this paper we are concerned with the nonhomogeneous heat conducting magnetohydrodynamic equations in $\mathbb{R}^2\times(0,T)$:
\begin{align}\label{mhd}
\begin{cases}
\rho_t + \divv(\rho u) = 0,\\
(\rho u)_t+\divv(\rho u\otimes u)-\mu\Delta u+\nabla P=H\cdot\nabla H,\\
c_{v}[(\rho \theta)_{t}+\divv(\rho u\theta)]-\kappa\Delta\theta =\frac{\mu}{2}|\nabla u+(\nabla u)^{tr}|^{2}+\nu(\curl H)^{2},\\
H_t-\nu\Delta H+u\cdot\nabla H-H\cdot\nabla u=0,\\
\divv u= \divv H=0.
\end{cases}
\end{align}
Here $\rho=\rho(x,t)$, $\theta=\theta(x,t)$, $u=(u^1,u^2)(x,t)$, $H=(H^1,H^2)(x,t)$, and $P=P(x,t)$ denote the density, the absolutely temperature, the velocity, the magnetic field, and the pressure, respectively. The positive constant $\mu$ is the viscosity coefficient of the fluid, $\nu>0$ is the magnetic diffusive coefficient, while
$c_v$ and $\kappa$ are the heat capacity and the ratio of the heat conductivity coefficient over the heat capacity, respectively.
$\curl H\triangleq\partial_1H^2-\partial_2H^1$.

The system \eqref{mhd} is supplemented with the initial condition
\begin{align}\label{a2}
(\rho, \rho u, \rho\theta, H)(x, 0)=(\rho_0, \rho_0u_0, \rho_0\theta_0, H_0)(x), \ \ x\in\Bbb R^2,
\end{align}
and the far field behavior
\begin{align}\label{n4}
(\rho, u, \theta, H)(x, t)\rightarrow (0, 0, 0, 0)\ \ {\rm as}\ |x|\rightarrow \infty, \ t>0.
\end{align}

Since the works of Lions \cite{L1996} and Choe-Kim \cite{CK2003},
where the global-in-time weak solutions and local strong solutions to the nonhomogeneous Navier-Stokes equations with vacuum (i.e., the initial density vanishes in some region) were obtained, respectively, there has been a considerable number of researches on the following nonhomogeneous magnetohydrodynamic equations in the presence of vacuum
\begin{align}\label{1.4}
\begin{cases}
\rho_t + \divv(\rho u) = 0,\\
(\rho u)_t+\divv(\rho u\otimes u)-\mu\Delta u+\nabla P=H\cdot\nabla H,\\
H_t-\nu\Delta H+u\cdot\nabla H-H\cdot\nabla u=0,\\
\divv u= \divv H=0.
\end{cases}
\end{align}
For a detailed derivation of the model \eqref{1.4}, we refer to \cite[Chapter 1]{GLL06}.
Under the compatibility condition
\begin{equation}\label{tan}
-\mu\Delta u+\nabla P_0-H_0\cdot \nabla H_0=\sqrt{\rho_0}g\ \ \text{for some}\ (P_0,g)\in H^1\times L^2,
\end{equation}
Chen-Tan-Wang \cite{tan} proved the local existence and uniqueness of strong solutions to the 3D Cauchy problem of \eqref{1.4}. At the same time, they obtained the global solution provided that the initial data satisfy some smallness condition. Later, with the help of a critical Sobolev inequality of logarithmic type involving the time, Huang and Wang \cite{hwjde1} derived the global strong solution in 2D bounded domains with general large initial data when the initial data satisfy \eqref{tan}. By virtue of spatial weighted estimates and the structural characteristic of \eqref{1.4},
L{\"u}-Xu-Zhong \cite{lvb2} established the global existence and uniqueness of strong solutions to the 2D Cauchy problem of \eqref{1.4} \textit{with vacuum at infinity}. Moreover, they also removed the compatibility condition \eqref{tan} by using time weighted techniques. Some important progress has been made about global strong solutions for the nonhomogeneous fluid equations with vacuum by many authors, please refer to \cite{HLL2021,DM2019,HW2014,Y19} and references therein. We apologize for not being able to list all the relevant references.

In contrast to \eqref{1.4}, the heat conducting model \eqref{mhd} is more in line with reality but the problem becomes challenging. It should be noted that \eqref{mhd} becomes the nonhomogeneous heat conducting Navier-Stokes equations when there is no electromagnetic field, we refer the reader to \cite[Chapter 2]{LK2016} for the detailed derivation of such system, and the mathematical results concerning the global existence of strong solutions to this model can refer for example to \cite{GL2021,Z22,ZT2015,Z2017}. Let's turn our attention to the system \eqref{mhd}. Wu \cite{W11} proved the local existence and uniqueness of strong solutions to the 3D initial boundary value problem of \eqref{mhd} provided that the initial data satisfy the compatibility condition
\begin{align}\label{1.6}
\begin{cases}
-\mu\Delta u_0+\nabla P_0-H_0\cdot\nabla H_0=\sqrt{\rho_0}g_1, \\
-\kappa\Delta\theta_0-\frac{\mu}{2}|\nabla u_0+(\nabla u_0)^{tr}|^2
-\nu(\curl H_0)^2=\sqrt{\rho_0}g_2,
\end{cases}
\end{align}
for some $P_0\in H^1$ and $g_1,g_2\in L^2$. This local well-posedness theory was very recently extended by Zhong \cite{zx4} to be a global one provided that $\big(\|\sqrt{\rho_0}u_0\|_{L^2}^2+\|H_0\|_{L^2}^2\big)
\big(\|\curl u_0\|_{L^2}^2+\|\curl H_0\|_{L^2}^2\big)$ is suitably small.
Such smallness condition is not needed to the 2D initial boundary value problem \cite{zx3} via Desjardins' interpolation inequality. Moreover, in \cite{zx3,zx4}, the author of this paper also proved that the velocity and the magnetic field converge exponentially to zero in $H^2$ and the gradient of the temperature converges algebraically to zero in $L^2$ as time goes to infinity, and there is no need to impose the compatibility condition \eqref{1.6} by applying time weighted techniques. We should point out that whether or not using the condition \eqref{1.6} may change with different problems. In \cite{zm2019}, \eqref{1.6} is required in order to ensure the boundedness of temperature when Zhu and Ou studied the global well-posedness of strong solutions for 3D initial boundary value problems with viscosity dependent density and temperature.
Meanwhile, to tackle the $L^\infty(0, T; L^2)$-norm of the gradient of the temperature, Zhong \cite{zx2} imposed the condition \eqref{1.6} and established global strong solution for large initial data to the 2D Cauchy problem of \eqref{mhd} with non-vacuum at infinity by a logarithmic interpolation inequality and delicate energy estimates.

Very recently, Chen and Zhong \cite{CZ21} showed the local existence and uniqueness of strong solutions to the problem \eqref{mhd}--\eqref{n4} with vacuum as far field density. However, the global well-posedness with general large initial data to \eqref{mhd}--\eqref{n4} \textit{with vacuum at infinity} is still open.
In fact, this is the main aim of the present paper.

Without loss of generality, we assume that the initial density $\rho_0$ satisfies
\begin{align}\label{oy3.7}
\int_{\mathbb{R}^2} \rho_0dx=1,
\end{align}
which implies that there exists a positive constant $N_0$ such that  \begin{align}\label{1.3}
\int_{B_{N_0}}\rho_0 dx\ge \frac12\int\rho_0dx=\frac12.
\end{align}
Here $B_{R}\triangleq\big\{x\in\mathbb{R}^2||x|<R\big\}$.

Our main result can be stated as follows.
\begin{theorem}\label{t1}
Let $\eta_0 $ be a positive constant and
\begin{align}\label{2.07}
\bar x\triangleq(e+|x|^2)^{\frac12}\ln^{1+\eta_0} (e+|x|^2).
\end{align}
For constants $q>2$ and $a>1$, in addition to \eqref{oy3.7},
assume that the initial data $(\rho_0\geq0, u_0, \theta_0\geq0, H_0)$ satisfies
\begin{align}\label{1.9}
\begin{cases}
\rho_0\bar x^a\in L^1 \cap H^1\cap W^{1,q},
\ H_0\bar x^{\frac{a}{2}}\in H^1, \\
\big(\sqrt{\rho_0}u_0,\sqrt{\rho_0}\theta_0\big)\in L^2,\
\big(\nabla u_0, \nabla \theta_0, \nabla H_0\big)\in H^1, \\
\divv u_0=\divv H_0=0,
\end{cases}
\end{align}
and the compatibility condition
\begin{align}\label{c.c}
\begin{cases}
-\mu\Delta u_0+\nabla P_0-H_0\cdot\nabla H_0=\sqrt{\rho_0}g_1, \\
\kappa\Delta\theta_0+\frac{\mu}{2}|\nabla u_0+(\nabla u_0)^{tr}|^2
+\nu(\curl H_0)^2=\sqrt{\rho_0}g_2,
\end{cases}
\end{align}
for some $P_0\in H^1(\mathbb{R}^2)$ and $g_1,g_2\in L^2(\mathbb{R}^2)$.
Then the problem  \eqref{mhd}--\eqref{n4} has a unique strong solution $(\rho\geq0,u,\theta\geq0,H)$ satisfying that, for any $0<T<\infty$,
\begin{align}\label{1.10}
\begin{cases}
\rho\bar x^a\in L^\infty(0,T;L^1\cap H^1\cap W^{1,q}),\\
\rho_t\in L^\infty(0,T;L^2\cap L^q),\\
\sqrt{\rho}u, \sqrt{\rho}\theta, \sqrt{\rho}u_t,
\sqrt{\rho}\theta_t, \nabla P\in L^\infty(0,T; L^2), \\
\nabla u, \nabla\theta, H\bar x^{\frac{a}{2}}\in\,L^\infty(0,T; H^1), \\
H, \nabla H, H_t, \nabla^2 H\in L^\infty(0,T;L^2), \\
\nabla u, \nabla\theta\in L^2(0,T;H^1)\cap L^2(0,T; W^{1,q})\cap  L^{\frac{q+1}{q}}(0,T; W^{1,q}), \\
\nabla P\in L^2(0,T;L^q)\cap L^{\frac{q+1}{q}}(0,T;L^q), \\
\nabla H, H_t, \nabla H\bar{x}^{\frac{a}{2}}\in L^2(0, T; H^1),\\
\sqrt{\rho}u_t, \sqrt{\rho}\theta_t, \nabla u_t,
\nabla\theta_t\in L^2(0, T;L^2),\\
\end{cases}
\end{align}
and
\begin{align}\label{l1.2}
\inf\limits_{0\le t\le T_0}\int_{B_{N_1}}\rho(x,t)dx\ge \frac14,
\end{align}
for some positive constant $N_1$ depending only on $\|\rho_0\|_{L^1}, \|\sqrt{\rho_0}u_0\|_{L^2}, \|H_0\|_{L^2}, N_0$, and $T$. Moreover, $(u,H)$ has the following decay rate, that is, for $t\geq1$,
\begin{align}\label{x5}
\|\nabla u(\cdot,t)\|_{L^2}+\|\nabla H(\cdot,t)\|_{L^2}\leq Ct^{-\frac12},
\end{align}
where $C$ depends only on $\mu$, $\nu$, $\|\rho_0\|_{L^\infty}$, $\|\sqrt{\rho_0}u_0\|_{L^2}$, $\|\nabla u_0\|_{L^2}$, and $\|H_0\|_{H^1}$.
\end{theorem}

\begin{remark}
It should be pointed out that the compatibility condition \eqref{c.c} is needed to obtain the $L^\infty(0, T; L^2)$-norm of $\sqrt{\rho}u_t$ and $\sqrt{\rho}\theta_t$, which is crucial in dealing with the $L^\infty(0, T; L^2)$-norm of the gradient of the temperature. It would be interesting to investigate whether such artificial condition could be removed such as in \cite{lvb2} via time weighted techniques.
\end{remark}

\begin{remark}
Due to the strong coupling between the velocity and the magnetic field,
we require the initial magnetic field to decay quickly at infinity. To our surprise, there is no need to impose such initial decay condition on the temperature although the temperature equation \eqref{mhd}$_3$ has a strong nonlinear term $|\nabla u+(\nabla u)^{tr}|^2$.
\end{remark}

\begin{remark}
Compared with \cite{zx2}, the decay rate \eqref{x5} is a new result. Moreover, we deduce from \eqref{2.4}, \eqref{3.4}, and \eqref{x5} that, for any $p\in[2,\infty)$ and $t\geq1$,
\begin{align*}
\|H\|_{L^p}\leq C\|H\|_{L^2}^{\frac2p}\|\nabla H\|_{L^2}^{\frac{p-2}{p}}
\leq Ct^{-\frac{p-2}{2p}}.
\end{align*}
We remark that it seems very hard to obtain the decay rate of the gradient of the temperature. The main difficulty lies in deriving time-independent spatial weighted estimate on the density (see \eqref{igj1-2}), which in turn effects a Hardy type estimate of the velocity (see \eqref{3.a2}).
\end{remark}

We now comment on the proof of Theorem \ref{t1}.
For the initial data satisfying \eqref{1.9} and \eqref{c.c}, the local existence and uniqueness of strong solutions to the problem \eqref{mhd}--\eqref{n4} has been established recently in \cite{CZ21} (see Lemma \ref{lem21}). Thus, one needs some global \textit{a priori} estimates on strong solutions to \eqref{mhd}--\eqref{n4} in suitable higher norms in order to extend the strong solution globally in time. It should be pointed out that the main difficulty here is the presence of \textit{vacuum at infinity} and the criticality of Sobolev's inequality in $\mathbb{R}^2$.
Technically, it seems difficult to bound the  $L^q(\mathbb{R}^2)$-norm of $u$
just in terms of $\|\sqrt{\rho}u\|_{L^{2}(\mathbb{R}^2)}$ and
$\|\nabla u\|_{L^{2}(\mathbb{R}^2)}$. Hence, the crucial techniques in \cite{W11} cannot be adapted because his arguments rely heavily on the fact that the $L^q$-norm of a function $u$ can be bounded by $\|\sqrt{\rho}u\|_{L^{2}}$ and $\|\nabla u\|_{L^{2}}$ for any $q\in[2,\infty)$ due to the absence of vacuum at infinity. Moreover, compared with \cite{lvb2}, some new difficulties arise due to the appearance of energy equation \eqref{mhd}$_3$ as well as the coupling of the velocity with the temperature. In fact, if we multiply \eqref{mhd}$_3$ by $\theta$ and integrate the resultant equality by parts over $\mathbb{R}^2$, then we have
\begin{align}\label{1.12}
\frac{c_v}{2}\frac{d}{dt}\int_{\mathbb{R}^2}\rho\theta^2dx
+\kappa\int_{\mathbb{R}^2}|\nabla \theta|^2dx
= \int_{\mathbb{R}^2}\Big[\frac{\mu}{2}|\nabla u+(\nabla u)^{tr}|^{2}+\nu(\curl H)^{2}\Big]\theta dx.
\end{align}
Since the $L^q(\mathbb{R}^2)$-norm of $\theta$ and spatial weighted estimates
on the gradients of the velocity and the magnetic field are unavailable, it is very hard to control the term on the right hand side of \eqref{1.12} directly.
To overcome this difficulty, motivated by \cite{L2021}, we establish a spatial weight estimate on the quadratic nonlinearity
$\frac{\mu}{2}|\nabla u+(\nabla u)^{tr}|^{2}+\nu(\curl H)^{2}$ (see \eqref{2.10}), which reveals that the gradients of the velocity and the magnetic field in a weighted $L^2$ space can be bounded by a weighted $L^1$-norm of $\rho\dot{\theta}$. This fact together with some estimates on
$(\rho,u,H)$ and a Hardy type estimate on $\theta$ implies that the right-hand side term of \eqref{1.12} can be controlled by the lower order norm of $\theta$ (see \eqref{3.8.3}). Fortunately, we can adopt similar strategies to tackle the $L^\infty(0,T;L^2(\mathbb{R}^2))$-norms of both $\nabla\theta$ and $\sqrt{\rho}\theta_t$ (see \eqref{3.77} and \eqref{x10}). Then, with the help of Gronwall's inequality, the $L_t^\infty L_x^2$-norms of $\sqrt{\rho}\theta$, $\nabla\theta$, and $\sqrt{\rho}\theta_t$
can be derived simultaneously by applying the compatibility condition \eqref{c.c}$_2$ and the \textit{a priori} estimates we
have obtained (see \eqref{3.8.11''}). Once with these estimates at hand,
the higher order bounds of the temperature can be shown by the standard $L^p$ theory of elliptic equations (see \eqref{3.8.14} and \eqref{3.8.16}). Finally, it is worth emphasizing that a Hardy type inequality (see \eqref{3.i2}) and Gagliardo-Nirenberg inequality (see \eqref{2.4} and \eqref{2.5}) are very useful for the analysis.

As a direct corollary of Theorem \ref{t1}, we have the following global existence result for 2D nonhomogeneous heat conducting Navier-Stokes equations with vacuum at infinity.
\begin{theorem}\label{co1}
Let $\eta_0 $ and $\bar x$ be as in \eqref{2.07}.
For constants $q>2$ and $a>1$, in addition to \eqref{oy3.7}, assume that the initial data $(\rho_0\geq0, u_0, \theta_0\geq0)$ satisfies
\begin{align*}
\rho_0\bar x^a\in L^1\cap H^1\cap W^{1,q},\
(\sqrt{\rho_0}u_0, \sqrt{\rho_0}\theta_0)\in L^2,\
(\nabla u_0,\nabla \theta_0)\in H^1, \
\divv u_0=0,
\end{align*}
and the compatibility condition
\begin{align*}
\begin{cases}
-\mu\Delta u_0+\nabla\widetilde{P}_0=\sqrt{\rho_0}\tilde{g}_1, \\
\kappa\Delta\theta_0+\frac{\mu}{2}|\nabla u_0+(\nabla u_0)^{tr}|^2
=\sqrt{\rho_0}\tilde{g}_2,
\end{cases}
\end{align*}
for some $\widetilde{P}_0\in H^1(\mathbb{R}^2)$ and $\tilde{g}_1,\tilde{g}_2\in L^2 (\mathbb{R}^2)$. Then, for any $0<T<\infty$, there exists a unique strong solution $(\rho,u,\theta)$ in $\mathbb{R}^2\times (0,T)$ to the 2D Cauchy problem of nonhomogeneous heat conducting Navier-Stokes equations (that is, \eqref{mhd}--\eqref{n4} with $H=0$) satisfying \eqref{1.10} and \eqref{x5} with $H=0$ and \eqref{l1.2}.
\end{theorem}

The rest of the paper is organized as follows. In Section \ref{sec2}, we collect some elementary facts and inequalities which will be needed in later analysis. Sections \ref{sec3} is devoted to the \textit{a priori} estimates which are needed to obtain the global existence of strong solutions. Finally, the main result Theorem \ref{t1} is proved in Section \ref{sec4}.

\section{Preliminaries}\label{sec2}

In this section, we will recall some  known facts and elementary
inequalities which will be used frequently later.
We begin with the local existence of strong solutions whose proof can be found in \cite{CZ21}.
\begin{lemma}\label{lem21}
Assume that $(\rho_0, u_0, \theta_0,H_0)$ satisfies \eqref{1.9} and \eqref{c.c}. Then there exists a small time $T>0$ and a unique strong solution $(\rho, u, \theta, H)$ to the problem \eqref{mhd}--\eqref{n4} in $\mathbb{R}^{2}\times(0,T]$ satisfying \eqref{1.10} and \eqref{l1.2}.
\end{lemma}

Next, the following well-known Gagliardo-Nirenberg inequality
(see \cite[Chapter II]{LSU1968}) will be used in the next section frequently.
\begin{lemma}
For $f\in H^1(\mathbb{R}^2)$ and $g\in L^r(\mathbb{R}^2)\cap W^{1, q}(\mathbb{R}^2)$ with $r\in(1,\infty)$ and $q\in(2,\infty)$, there exists a positive constant $C$ such that
\begin{align}\label{2.4}
&\|f\|^p_{L^p}\leq C\|f\|^{2}_{L^2}\|\nabla f\|^{p-2}_{L^2}, \ \forall p\in[2,\infty),\\ \label{2.5}
&\|g\|_{L^\infty}\leq
C\|g\|^{\frac{r(q-2)}{2q+r(q-2)}}_{L^r} \|\nabla g\|^{\frac{2q}{2q+r(q-2)}}_{L^q}.
\end{align}
\end{lemma}

Next, for $\Omega=\mathbb{R}^2$ or $\Omega=B_R$, the following weighted $L^m$-bounds for elements of the Hilbert space $\tilde{D}^{1,2}(\Omega)\triangleq\{v\in H^1_{\rm loc}(\Omega)|\nabla v\in L^2(\Omega)\}$ can be found in \cite[Theorem B.1]{L1996}.
\begin{lemma}\label{1leo}
For $m\in [2,\infty)$ and $\theta\in (1+\frac{m}{2},\infty),$ there exists a positive constant $C$ such that for either $\Omega=\mathbb{R}^2$ or $\Omega=B_R$ with $R\ge 1$ and for any $v\in \widetilde{D}^{1,2}(\Omega)$,
\begin{align}\label{3h}
\left(\int_{\Omega} \frac{|v|^m}{(e+|x|^2)\ln^{\theta}(e+|x|^2)}dx  \right)^{\frac1m}\le C\|v\|_{L^2(B_1)}+C\|\nabla v\|_{L^2 (\Omega)}.
\end{align}
\end{lemma}

A useful consequence of Lemma \ref{1leo} is the following crucial weighted  bounds (see \cite[Lemma 2.4]{L2015}) for elements of $\widetilde{D}^{1,2}(\Omega)$.
\begin{lemma}\label{lemma2.3}
Let $\bar x$ and $\eta_0$ be as in \eqref{2.07} and $\Omega$ be as in Lemma \ref{1leo}. Assume that $\rho\in L^\infty(\Omega)$ is a non-negative function such that
\begin{align}\label{2.i2}
\int_{B_{N_1}}\rho dx\ge M_1, \ \|\rho\|_{L^\infty(\Omega)}\le M_2,
\end{align}
for positive constants $M_1, M_2$, and $N_1\ge 1$ with $B_{N_1}\subset\Omega.$ Then, for $\varepsilon,\eta>0$, there is a positive constant $C$ depending only on $\varepsilon,\eta, M_1,M_2,
N_1,$ and $\eta_0$ such that, for $v\in \widetilde{D}^{1,2}(\Omega)$ with $\sqrt{\rho}v\in L^2(\Omega)$,
\begin{align}\label{3.i2}
\|v\bar x^{-\eta}\|_{L^{\frac{2+\varepsilon}{\tilde{\eta}}}(\Omega)}
&\le C\big(\|\sqrt{\rho}v\|_{L^2(\Omega)}+\|\nabla v\|_{L^2(\Omega)}\big)
\end{align}
with $\tilde{\eta}=\min\{1,\eta\}$.
\end{lemma}

Let $\mathcal{H}^{1}(\mathbb{R}^2)$ and $BMO(\mathbb{R}^2)$ stand for the usual Hardy and $BMO$ spaces (see \cite[Chapter IV]{S1993}). Then the following well-known facts play a key role in the proof of Lemma \ref{l3.0} in the next section.
\begin{lemma}\label{lem27}
(a) There is a positive constant $C$ such that
\begin{equation*}
\|E\cdot B\|_{\mathcal{H}^{1}(\mathbb{R}^2)}
\leq C\|E\|_{L^{2}(\mathbb{R}^2)}\|B\|_{L^{2}(\mathbb{R}^2)},
\end{equation*} for all $E\in L^{2}(\mathbb{R}^2)$ and $ B\in L^{2}(\mathbb{R}^2)$ satisfying
\begin{equation*}
\divv E=0,\ \nabla^{\bot}\cdot B=0\ \ \text{in}\ \ \mathcal{D}'(\mathbb{R}^2).
\end{equation*}
(b) There is a positive constant $C$ such that
\begin{equation}\label{lem1}
\| v \|_{BMO(\mathbb{R}^2)}\leq C\|\nabla v \|_{L^{2}(\mathbb{R}^2)},
\end{equation} for all $  v \in D^1(\mathbb{R}^2)$.
\end{lemma}
{\it Proof.}
(a) For the detailed proof, please see \cite[Theorem II.1]{CLMS1993}.

(b) It follows  from the Poincar{\'e} inequality that for any ball $B\subset\mathbb{R}^2$
\begin{equation*}
\frac{1}{|B|}\int_{B}\left| v(x) - \frac{1}{|B|}\int_Bv(y)dy\right|dx\leq C\left(\int_{B}|\nabla v |^2dx\right)^{\frac12},
\end{equation*}
which directly gives \eqref{lem1}.\hfill $\Box$

Finally, we have the following spatial weighted estimate which is very important in dealing with the \textit{a priori} estimates of the temperature.
\begin{lemma}\label{lemma2.5}
Let $(\rho,u,\theta,H)$ be the solution to the problem \eqref{mhd}--\eqref{n4}, then it holds that, for any $b_1>0$,
\begin{align}\label{2.10}
\int_{\mathbb{R}^2} \Big[\frac{\mu}{2}|\nabla u+(\nabla u)^{tr}|^{2}+\nu(\curl H)^{2}\Big]|x|^{b_1}dx
\leq \int_{\mathbb{R}^2}[c_v(\rho\theta_t+\rho u\cdot\nabla\theta)]|x|^{b_1}dx.
\end{align}
\textit{Proof.} Applying standard maximum principle (see \cite[p. 43]{F2004}) to \eqref{mhd}$_3$ along with $\theta_0\geq0$ shows that
\begin{equation}\label{x7}
\inf_{\mathbb{R}^2\times[0,T]}\theta(x,t)\geq0.
\end{equation}
For $b_1>0$, direct calculation gives that, for $i=1,2$,
\begin{align*}
\partial_{x_ix_i}|x|^{b_1}=b_1|x|^{b_1-2}+b_1(b_1-2)x_i^2|x|^{b_1-4}.
\end{align*}
This implies that
\begin{align}\label{x8}
\Delta|x|^{b_1} & =\partial_{x_1x_1}|x|^{b_1}+\partial_{x_2x_2}|x|^{b_1} \notag \\
& =2b_1|x|^{b_1-2}+b_1(b_1-2)|x|^{b_1-4}(x_1^2+x_2^2) \notag \\
& =b_1^2|x|^{b_1-2}.
\end{align}
Multiplying \eqref{mhd}$_3$ by $|x|^{b_1}$ and integrating the resultant equality over $\mathbb{R}^2$, we obtain from \eqref{mhd}$_1$ that
\begin{align}\label{x9}
\int_{\mathbb{R}^2}[c_v(\rho\theta_t+\rho u\cdot\nabla\theta)]|x|^{b_1}dx
=\int_{\mathbb{R}^2} \Big[\frac{\mu}{2}|\nabla u+(\nabla u)^{tr}|^{2}+\nu(\curl H)^{2}\Big]|x|^{b_1}dx+\kappa\int_{\mathbb{R}^2}\Delta\theta|x|^{b_1}dx.
\end{align}
Integration by parts together with \eqref{x7} and \eqref{x8} yields that
\begin{align*}
\kappa\int_{\mathbb{R}^2}\Delta\theta|x|^{b_1}dx
=\kappa\int_{\mathbb{R}^2}\theta\Delta|x|^{b_1}dx
=\kappa b_1^2\int_{\mathbb{R}^2}\theta|x|^{b_1-2}dx\geq0,
\end{align*}
which combined with \eqref{x9} implies \eqref{2.10}.
\end{lemma}

\section{\textit{A priori} estimates}\label{sec3}

In this section, we will establish some necessary \textit{a priori} bounds for strong solutions $(\rho,u,\theta,H)$ to the problem \eqref{mhd}--\eqref{n4} to extend the local strong solution. Thus, let $T>0$ be a fixed time and $(\rho,u,\theta,H)$  be the strong solution to \eqref{mhd}--\eqref{n4} on $\mathbb{R}^2\times(0,T]$ with initial data $(\rho_0,u_0,\theta_0,H_0)$ satisfying \eqref{1.9} and \eqref{c.c}. In what follows, for simplicity, we write
\begin{equation*}
\int\cdot dx=\int_{\mathbb{R}^2}\cdot dx.
\end{equation*}
Moreover, we sometimes use $C(\alpha)$ to emphasize the dependence on $\alpha$.

We begin with the following elementary estimate for the solution.
\begin{lemma} \label{l3.00}
It holds that
\begin{align}\label{3.1}
& \sup_{0\le t\le T}\big(\|\rho\|_{L^1\cap L^\infty}+\|\sqrt{\rho}u\|^2_{L^{2}}+\|H\|_{L^2}^2\big)
+\int_{0}^{T}\big(\mu\|\nabla u\|_{L^2}^2+\nu\|\nabla H\|_{L^2}^2\big)dt \notag \\
& \le \|\rho_0\|_{L^{1}\cap L^\infty}
+\|\sqrt{\rho_0}u_0\|^2_{L^2}+\|H_0\|_{L^2}^2.
\end{align}
\end{lemma}
\textit{Proof.}
We deduce from \eqref{mhd}$_1$ and \eqref{mhd}$_5$ that the density satisfies a transport equation, thus we have
\begin{equation}\label{3.2}
\sup_{0\le t\le T}\|\rho\|_{L^p}
\leq\|\rho_0\|_{L^p},\ \forall\ 1\leq p\leq \infty.
\end{equation}
Multiplying \eqref{mhd}$_2$ by $u$ and \eqref{mhd}$_4$ by $H$, respectively, we get after integrating by parts that
\begin{align}\label{3.3}
\frac{d}{dt}\big(\|\sqrt{\rho}u\|_{L^2}^2+\|H\|_{L^2}^2\big)
+2\big(\mu\|\nabla u\|_{L^2}^2+\nu\|\nabla H\|_{L^2}^2\big)
=0.
\end{align}
Integrating \eqref{3.3} over $[0,T]$ leads to
\begin{align}\label{3.4}
\sup_{0\le t\le T} \big(\|\sqrt{\rho}u\|^2_{L^2}+\|H\|_{L^2}^2\big)
+\int_{0}^{T}\big(\mu\|\nabla u\|_{L^2}^2+\nu\|\nabla H\|_{L^2}^2\big)dt
\le \|\sqrt{\rho_0}u_0\|^2_{L^2}+\|H_0\|_{L^2}^2.
\end{align}
This together with \eqref{3.2} yields the desired \eqref{3.1}.
\hfill $\Box$

Next, the following lemma concerns the key uniformly-in-time estimate on the $L^\infty(0,T;L^2)$-norm of the gradients of the velocity and the magnetic field.
\begin{lemma}\label{l3.0}
There exists a positive constant $C$ depending only on $\mu$, $\nu$, $\|\rho_0\|_{L^\infty}$, $\|\sqrt{\rho_0}u_0\|_{L^2}$, $\|\nabla u_0\|_{L^2}$, and $\|H_0\|_{H^1}$ such that
\begin{align}\label{3.5}
\sup_{0\le t\le T}\big(\|\nabla u\|_{L^2}^2+\|\nabla H\|_{L^2}^2\big)
+\int_{0}^{T}\big(\|\sqrt{\rho}\dot{u}\|_{L^2}^2+\|\nabla^2u\|_{L^2}^2
+\||H||\nabla H|\|_{L^2}^2+\|\nabla^2H\|_{L^2}^2\big)dt
\le C.
\end{align}
Here $\dot{v}\triangleq\partial_{t}v+u\cdot\nabla v$. Moreover, one has
\begin{align}\label{x4}
\sup_{0\le t\le T}\big[t\big(\|\nabla u\|_{L^2}^2+\|\nabla H\|_{L^2}^2\big)\big]\le C.
\end{align}
\end{lemma}
\textit{Proof.}
1. Multiplying \eqref{mhd}$_2$ by $\dot{u}$ and integrating the resulting equality over $\mathbb{R}^2$ lead  to
\begin{align}\label{3.6}
\int\rho|\dot{u}|^{2}dx
&=\int \mu\Delta u\cdot\dot{u}dx-\int\nabla P\cdot\dot{u} dx
-\frac12\int\nabla|H|^2\cdot\dot udx+\int H\cdot\nabla H\cdot\dot udx
\notag \\ &\triangleq I_{1}+I_{2}+I_{3}+I_{4}.
\end{align}
It follows from integration by parts and \eqref{2.4} that
\begin{align}\label{3.7}
I_{1} & =\int\mu\Delta u\cdot(u_t+u\cdot\nabla u)dx \notag \\
& = -\frac{\mu}{2}\frac{d}{dt}\|\nabla u\|_{L^2}^2
-\mu\int\partial_{i}u^{j}\partial_{i}(u^{k}\partial_{k}u^{j})dx  \notag \\
& \leq-\frac{\mu}{2}\frac{d}{dt}\|\nabla u\|_{L^2}^2
+C\|\nabla u\|_{L^3}^3 \notag \\
& \leq-\frac{\mu}{2}\frac{d}{dt}\|\nabla u\|_{L^2}^2
+C\|\nabla u\|_{L^2}^2\|\nabla^2u\|_{L^2}.
\end{align}
Here and in what follows, we use the Einstein convention that the repeated indices denote the summation.
Integration by parts together with \eqref{mhd}$_5$ gives rise to
\begin{align*}
I_{2}=-\int\nabla P\cdot(u_t+u\cdot\nabla u)dx
= \int P\partial_{j}u^{i}\partial_{i}u^{j}dx
\leq C\|P\|_{BMO}\|\partial_{j}u^{i}\partial_{i}u^{j}\|_{\mathcal{H}^1},
\end{align*}
where one has used the duality of $\mathcal{H}^1$ and $BMO$ (see \cite[Chapter IV]{S1993}) in the last inequality. Since $\divv(\partial_j u )=\partial_j\divv u=0$ and $\nabla^{\bot}\cdot(\nabla u^{j})=0$, then Lemma \ref{lem27} yields that
\begin{equation}
I_{2}\leq C\|P\|_{BMO}\|\partial_{j}u^{i}\partial_{i}u^{j}\|_{\mathcal{H}^1}
\leq C\|\nabla P\|_{L^2}\|\nabla u\|_{L^2}^{2}.
\end{equation}
For the term $I_3$, integration by parts together with \eqref{mhd}$_5$ and \eqref{2.4} leads to
\begin{align}
I_3=\frac12\int|H|^2\partial_iu^j\partial_ju^idx
\leq C\|H\|_{L^6}^6+C\|\nabla u\|_{L^3}^3
\leq C\|H\|_{L^2}^2\|\nabla H\|_{L^2}^4
+C\|\nabla u\|_{L^2}^{2}\|\nabla^2u\|_{L^2}.
\end{align}
Integration by parts, we infer from \eqref{mhd}$_4$, \eqref{mhd}$_5$, H{\"o}lder's inequality, \eqref{2.4}, and Young's inequality that, for $\delta>0$,
\begin{align}\label{3.10}
I_4
& = \int H\cdot\nabla H\cdot u_tdx+\int H\cdot\nabla H\cdot (u\cdot \nabla u)dx
\notag \\
& = -\frac{d}{dt}\int H\cdot\nabla u\cdot Hdx+\int H_t\cdot\nabla u\cdot Hdx+\int H\cdot\nabla u\cdot  H_tdx \notag \\
& \quad- \int H^i\partial_i u^j\partial_j u^k H^kdx- \int H^i u^j\partial_i\partial_j u^k H^kdx
\notag \\
& =-\frac{d}{dt}\int H\cdot\nabla u \cdot  Hdx+\int (\nu\Delta H-u\cdot\nabla H+H\cdot\nabla u)\cdot\nabla u\cdot  Hdx \notag \\
& \quad+\int H\cdot\nabla u\cdot (\nu\Delta H-u\cdot\nabla H+H\cdot\nabla u)dx- \int H^i\partial_i u^j\partial_j u^k H^kdx \notag \\
& \quad+\int u^j\partial_jH^i\partial_i u^k H^kdx +\int H^i \partial_i u^k u^j\partial_jH^kdx
\notag \\
& =-\frac{d}{dt}\int H\cdot\nabla u \cdot Hdx+\nu\int\Delta H\cdot\nabla u\cdot Hdx
+\nu\int H\cdot\nabla u\cdot\Delta Hdx+\int H\cdot\nabla u\cdot H\cdot\nabla udx
\notag \\
& \le -\frac{d}{dt}\int H\cdot\nabla u\cdot Hdx
+2\nu\|\Delta H\|_{L^2}\|H\|_{L^6}\|\nabla u\|_{L^3}
+\|H\|_{L^6}^2\|\nabla u\|_{L^3}^2 \notag \\
& \le -\frac{d}{dt}\int H\cdot\nabla u\cdot Hdx
+\frac\delta2\|\Delta H\|_{L^2}^2
+C\|H\|_{L^2}^{\frac23}\|\nabla H\|_{L^2}^{\frac43}
\|\nabla u\|_{L^2}^{\frac43}\|\nabla^2u\|_{L^2}^{\frac23} \notag \\
& \le -\frac{d}{dt}\int H\cdot\nabla u\cdot Hdx+\frac\delta2\|\Delta H\|_{L^2}^2
+C\|H\|_{L^2}^2\|\nabla H\|_{L^2}^4+C\|\nabla u\|_{L^2}^{2}\|\nabla^2u\|_{L^2}.
\end{align}
Hence, inserting \eqref{3.7}--\eqref{3.10} into \eqref{3.6} and using \eqref{3.4}, we arrive at
\begin{align}\label{3.11}
B'(t)+\|\sqrt{\rho}\dot{u}\|_{L^2}^2
\leq \delta\|\Delta H\|_{L^2}^2+C\|\nabla H\|_{L^2}^4
+C\big(\|\nabla^2u\|_{L^2}+\|\nabla P\|_{L^2}\big)\|\nabla u\|_{L^2}^2.
\end{align}
where
$$B(t)\triangleq\mu\|\nabla u\|_{L^2}^2+2\int H\cdot\nabla u \cdot Hdx$$
satisfies
\begin{align}\label{3.12}
\frac{\mu}{2}\|\nabla u\|_{L^2}^2-C_1\|\nabla H\|_{L^2}^2
\leq B(t)\leq\frac{3\mu}{2}\|\nabla u\|_{L^2}^2+C_1\|\nabla H\|_{L^2}^2
\end{align}
for some positive constant $C_1$ depending only on $\mu$
due to the following fact
\begin{align*}
\left|2\int H\cdot\nabla u \cdot Hdx\right|
& \leq 2\|\nabla u\|_{L^2}\|H\|_{L^4}^2 \\
& \leq \frac{\mu}{2}\|\nabla u\|_{L^2}^2+\frac{2}{\mu}\|H\|_{L^4}^4 \\
& \leq \frac{\mu}{2}\|\nabla u\|_{L^2}^2+C(\mu)\|H\|_{L^2}^2\|\nabla H\|_{L^2}^2 \\
& \leq \frac{\mu}{2}\|\nabla u\|_{L^2}^2+C_1\|\nabla H\|_{L^2}^2.
\end{align*}

2. Multiplying \eqref{mhd}$_4$ by $\Delta H$ and integrating the resultant equations by parts over $\mathbb{R}^2$,
it follows from H\"older's inequality, \eqref{2.4}, \eqref{3.4}, and Young's inequality that
\begin{align}\label{3.13}
& \frac{d}{dt}\int|\nabla H|^2 dx+2\nu\int|\Delta H|^2dx \notag \\
& \leq C\int |\nabla u||\nabla H|^2dx+C\int|\nabla u||H||\Delta H|dx \notag \\
& \leq C\|\nabla u\|_{L^3}\|\nabla H\|_{L^2}\|\nabla H\|_{L^6}
+C\|\nabla u\|_{L^3}\|H\|_{L^6}\|\Delta H\|_{L^2}\notag \\
& \leq C\|\nabla u\|_{L^2}^{\frac23}\|\nabla^2u\|_{L^2}^{\frac13}
\|\nabla H\|_{L^2}^{\frac43}\|\nabla^2H\|_{L^2}^{\frac23}
+C\|\nabla u\|_{L^2}^{\frac23}\|\nabla^2u\|_{L^2}^{\frac13}
\|H\|_{L^2}^{\frac13}\|\nabla H\|_{L^2}^{\frac23}\|\Delta H\|_{L^2}\notag \\
& \leq C\|\nabla u\|_{L^2}^2\|\nabla^2u\|_{L^2}
+C\|\nabla H\|_{L^2}^2\|\nabla^2H\|_{L^2}
+C\|\nabla H\|_{L^2}\|\Delta H\|_{L^2}^{\frac32}\notag \\
& \leq C\|\nabla u\|_{L^2}^2\|\nabla^2u\|_{L^2}
+C\|\nabla H\|_{L^2}^4+\nu\|\Delta H\|_{L^2}^2,
\end{align}
where we have used
\begin{align}\label{x}
\|\nabla^2H\|_{L^2}\leq C\|\Delta H\|_{L^2},
\end{align}
due to the standard $L^2$-estimate of elliptic equations.
Thus, adding \eqref{3.13} multiplied by $2C_1$ to \eqref{3.11} and choosing  $\delta=C_1\nu$, we get that
\begin{align}\label{3.14}
& \frac{d}{dt}\big(B(t)
+2C_1\|\nabla H\|_{L^2}^2\big)+\|\sqrt{\rho}\dot{u}\|_{L^2}^2
+C_1\nu\|\Delta H\|_{L^2}^2 \notag \\
& \le C\|\nabla H\|_{L^2}^4
+C\big(\|\nabla^2u\|_{L^2}+\|\nabla P\|_{L^2}\big)\|\nabla u\|_{L^2}^2.
\end{align}

3. Noting that $(u,P)$ satisfies the Stokes system
\begin{align}\label{3.15}
\begin{cases}
-\mu\Delta u+\nabla P= -\rho\dot{u}+H\cdot\nabla H,\ & x\in \mathbb{R}^2,\\
\divv u=0, \ & x\in  \mathbb{R}^2,\\
u(x)=0,\  & |x|\rightarrow\infty.
\end{cases}
\end{align}
Applying the standard $L^p$-estimate to \eqref{3.15} yields that, for any $p\in [2,\infty)$,
\begin{align}\label{3.16}
\|\nabla^2u\|_{L^p}+\|\nabla P\|_{L^p}
\le C\|\rho\dot{u}\|_{L^p}+C\||H||\nabla H|\|_{L^p}.
\end{align}
Then we obtain from \eqref{3.16} with $p=2$ and \eqref{3.2} that
\begin{align}\label{3.17}
\|\nabla^2u\|_{L^2}+\|\nabla P\|_{L^2}
\leq C\|\rho\dot{u}\|_{L^2}+C\||H||\nabla H|\|_{L^2}
\leq C\|\sqrt{\rho}\dot{u}\|_{L^2}+C\||H||\nabla H|\|_{L^2}.
\end{align}
Putting \eqref{3.17} into \eqref{3.14} and applying Cauchy-Schwarz inequality, we have
\begin{align}\label{3.18}
& \frac{d}{dt}\big(B(t)+2C_1\|\nabla H\|_{L^2}^2\big)
+\frac12\|\sqrt{\rho}\dot{u}\|_{L^2}^2
+C_1\nu\|\Delta H\|_{L^2}^2 \notag \\
& \le C\big(\|\nabla H\|_{L^2}^2+\|\nabla u\|_{L^2}^2\big)\big(\|\nabla H\|_{L^2}^2+\|\nabla u\|_{L^2}^2\big)+C\||H||\nabla H|\|_{L^2}^2.
\end{align}
Multiplying\eqref{mhd}$_4$ by $|H|^2H$ and integrating the resultant equality over $\mathbb{R}^2$, we obtain from \eqref{2.4} that
\begin{align}\label{x2}
\frac14\frac{d}{dt}\|H\|^4_{L^4}+\||\nabla H| |H|\|_{L^2}^2
+\frac12\|\nabla|H|^2\|_{L^2}^2
& \leq C\|\nabla u\|_{L^2} \||H|^2\|_{L^4}^2 \notag \\
&\leq C\|\nabla u\|_{L^2} \||H|^2\|_{L^2}\|\nabla|H|^2\|_{L^2} \notag \\
&\leq \frac14\|\nabla|H|^2\|_{L^2}^2
+C\|\nabla u\|_{L^2}^2\|H\|_{L^4}^4,
\end{align}
which together with Gronwall's inequality and \eqref{3.4} implies that
\begin{align}\label{3.19}
& \sup_{0\leq t\leq T}\|H\|_{L^4}^4
+\int_{0}^{T}\||H||\nabla H|\|_{L^2}^2dt \leq C.
\end{align}
Hence, \eqref{3.5} follows from \eqref{3.18}, Gronwall's inequality, \eqref{3.12}, \eqref{3.17}, \eqref{x}, \eqref{3.4}, and \eqref{3.19}.

4. Multiplying \eqref{x2} by $t$, we then obtain from \eqref{2.4} and \eqref{3.4} that
\begin{align*}
\frac{d}{dt}\big(t\|H\|^4_{L^4}\big)+t\||\nabla H| |H|\|_{L^2}^2
&\leq C\|\nabla u\|_{L^2}^2\big(t\|H\|^4_{L^4}\big)+C\|H\|_{L^4}^4 \\
&\leq C\|\nabla u\|_{L^2}^2\big(t\|H\|^4_{L^4}\big)+C\|H\|_{L^2}^2\|\nabla H\|_{L^2}^2\\
&\leq C\|\nabla u\|_{L^2}^2\big(t\|H\|^4_{L^4}\big)+C\|\nabla H\|_{L^2}^2,
\end{align*}
which combined with Gronwall's inequality and \eqref{3.4} yields that
\begin{align}\label{x3}
& \sup_{0\leq t\leq T}\big(t\|H\|_{L^4}^4\big)
+\int_{0}^{T}t\||H||\nabla H|\|_{L^2}^2dt \leq C.
\end{align}
Multiplying \eqref{3.18} by $t$ and using \eqref{3.12}, we have
\begin{align*}
& \frac{d}{dt}\big[t\big(B(t)+2C_1\|\nabla H\|_{L^2}^2\big)\big]
+\frac12t\|\sqrt{\rho}\dot{u}\|_{L^2}^2
+C_1\nu t\|\Delta H\|_{L^2}^2 \notag \\
& \le C\big(\|\nabla H\|_{L^2}^2+\|\nabla u\|_{L^2}^2\big)\big[t\big(\|\nabla H\|_{L^2}^2+\|\nabla u\|_{L^2}^2\big)\big]
+Ct\||H||\nabla H|\|_{L^2}^2+C\big(\|\nabla u\|_{L^2}^2+\|\nabla H\|_{L^2}^2\big).
\end{align*}
This along with Gronwall's inequality, \eqref{3.12}, \eqref{x3}, and \eqref{3.4} implies \eqref{x4}.
\hfill $\Box$

The following spatial weighted estimates on the density and the magnetic field play an important role in bounding the higher order derivatives of the solution.
\begin{lemma} \label{l3.01}
There exists a positive constant $C$ depending on $T$ such that
\begin{align}\label{igj1}
\sup_{0\le t\le T} \big(\|\rho\bar x^a\|_{L^1}+\|H\bar x^{\frac{a}{2}}\|_{L^2}^2\big)+\int_{0}^{T}\|\nabla H \bar x^{\frac{a}{2}}\|_{L^2}^2 dt\le C(T).
\end{align}
\end{lemma}
\textit{Proof.} 1. For $N>1$, let $\varphi_N\in C^\infty_0(B_N)$ satisfy
\begin{align}\label{vp1}
0\le \varphi_N \le 1, \ \varphi_N(x)=1,\ \mbox{if}\ |x|\le \frac{N}{2},\ \text{and} \ |\nabla \varphi_N|\le 3N^{-1}.
\end{align}
It follows from \eqref{mhd}$_1$, \eqref{3.2}, and \eqref{3.4} that
\begin{align}\label{oo0}
\frac{d}{dt}\int \rho\varphi_{2N_0}dx
& =\int\rho u\cdot\nabla\varphi_{2N_0}dx
\ge -CN_0^{-1}\|\rho\|_{L^1}^{\frac12}
\|\sqrt{\rho}u\|_{L^2}\ge -\widetilde{C}N_0^{-1}
\end{align}
for some positive constant $\widetilde{C}$ depending only on
$\|\rho_0\|_{L^1}$, $\|\sqrt{\rho_0}u_0\|_{L^2}$, and $\|H_0\|_{L^2}$.
Integrating \eqref{oo0} with respect to the time and choosing
$N= N_1\triangleq2N_0+4\widetilde{C}T$, we obtain after using \eqref{1.3} that
\begin{align}\label{p1}
\inf\limits_{0\le t\le T}\int_{B_{N_1}} \rho dx
&\ge \inf\limits_{0\le t\le T}\int \rho \varphi_{N_1} dx
\ge \int \rho_0 \varphi_{N_1} dx-\widetilde{C} N_1^{-1}T
\ge \int_{B_{N_0}} \rho_0dx-\frac{\widetilde{C} T}{2N_0+4\widetilde{C} T}
\ge \frac14.
\end{align}
The combination of \eqref{p1}, \eqref{3.1}, and \eqref{3.i2} implies that, for $\varepsilon,\eta>0$ and $v\in \widetilde{D}^{1,2}(B_R)$ with $\sqrt{\rho}v\in L^2(B_R)$,
\begin{align}\label{3.v2}
\|v\bar x^{-\eta}\|_{L^{\frac{2+\varepsilon}{\tilde{\eta}}}}^2
&\le C(\varepsilon,\eta)\big(\|\sqrt{\rho}v\|_{L^2}^2
+\|\nabla v\|_{L^2}^2\big),
\end{align}
where $\tilde{\eta}\triangleq\min\{1,\eta\}.$

2. Noting that for any $s>0$, it holds that
\begin{align}\label{z10.6}
|\nabla\bar{x}|\leq C(\eta_0)\ln^{1+\eta_0}(e+|x|^2)
\leq C(\eta_0)\bar{x}^{s}.
\end{align}
Multiplying $\eqref{mhd}_{1}$ by $\bar x^a$ and integrating by parts, we then obtain from H{\"o}lder's inequality, \eqref{3.v2}, \eqref{z10.6}, and \eqref{3.1} that
\begin{align*}
\frac{d}{dt}\|\rho\bar x^a\|_{L^1}
& = \int\rho(u\cdot\nabla)\bar xa\bar x^{a-1}dx \\
& \le C \int \rho|u|\bar x^{a-1+\frac{4}{8+a}}dx \notag \\
& \le C\|\rho\bar x^{a-1+\frac{8}{8+a}}\|_{L^\frac{8+a}{7+a}}
\|u\bar x^{-\frac{4}{8+a}}\|_{L^{8+a}} \notag \\
& \le C\|\rho\|_{L^\infty}^{\frac{1}{8+a}}\|\rho\bar x^a\|_{L^1}^{\frac{7+a}{8+a}}
\big(\|\sqrt{\rho}u\|_{L^2} +\|\nabla u\|_{L^2}\big) \notag \\
& \le C\|\rho\bar x^a\|_{L^1}+C.
\end{align*}
This combined with Gronwall's inequality and \eqref{1.9} leads to
\begin{align}\label{igj1-2}
\sup_{0\le t\le T}\|\rho\bar x^a\|_{L^1}\le C(T).
\end{align}
It follows from H{\"o}lder's inequality, \eqref{3.2}, \eqref{3.v2}, and \eqref{igj1-2} that, for any $\varepsilon,\eta>0$ and $v\in \widetilde{D}^{1,2}(B_R)$ with $\sqrt{\rho}v\in L^2(B_R)$,
\begin{align}\label{local1}
\|\rho^\eta v\|_{L^{\frac{2+\varepsilon}{\tilde{\eta}}}}
& \le C\|\rho^\eta \bar x^{\frac{3\tilde{\eta} a}{4(2+\varepsilon)}} \|_{L^{ \frac{4(2+\varepsilon)}{3\tilde{\eta}}}} \|v\bar x^{-\frac{3\tilde{\eta} a}{4(2+\varepsilon)}} \|_{L^{ \frac{4(2+\varepsilon)}{\tilde{\eta}}}} \notag \\
& \le C\left(\int \rho^{\frac{4(2+\varepsilon)\eta}{3\tilde{\eta}}-1}\rho\bar x^a dx\right)^{ \frac{3\tilde{\eta}}{4(2+\varepsilon)}} \|v\bar x^{-\frac{3\tilde{\eta} a}{4(2+\varepsilon)}} \|_{L^{ \frac{4(2+\varepsilon)}{\tilde{\eta}}}} \notag \\
& \le C\|\rho\|_{L^\infty}^{\frac{4(2+\varepsilon)\eta-3\tilde{\eta}}{4(2+\varepsilon)}}
\|\rho\bar x^a\|_{L^1}^{\frac{3\tilde{\eta}}{4(2+\varepsilon)}}
\big(\|\sqrt{\rho}v\|_{L^2}+\|\nabla v\|_{L^2}\big) \notag \\
& \le C\|\sqrt{\rho}v\|_{L^2}+C\|\nabla v\|_{L^2},
\end{align}
where $\tilde{\eta}=\min\{1,\eta\}$. In particular, this together with \eqref{3.v2}, \eqref{3.4}, and \eqref{3.5} implies that
\begin{align}\label{3.a2}
&\|\rho^\eta u\|_{L^{\frac{2+\varepsilon}{\tilde{\eta}}}}
+\|u\bar x^{-\eta}\|_{L^{\frac{2+\varepsilon}{\tilde{\eta}}}}
\le C\big(\|\sqrt{\rho}u\|_{L^2}+\|\nabla u\|_{L^2}\big)\le C.
\end{align}

3. Multiplying \eqref{mhd}$_4$ by $H\bar{x}^a$ and integrating by parts yield
\begin{align}\label{lv4.1}
\frac{1}{2}\frac{d}{dt}\|H\bar{x}^{\frac{a}{2}}\|_{L^2}^2+\nu \|\nabla H \bar{x}^{\frac{a}{2}}\|_{L^2}^2
& =\frac{\nu}{2}\int |H|^2\Delta\bar{x}^adx+\int (H\cdot\nabla)u\cdot H\bar{x}^adx +\frac12\int |H|^2u\cdot\nabla\bar{x}^adx \notag \\
& \triangleq \bar{I}_1+\bar{I}_2+\bar{I}_3.
\end{align}
Direct calculations lead to
\begin{equation}\label{10.4}
|\bar{I}_1|\leq C\int|H|^2\bar{x}^a\bar{x}^{-1}dx
\leq C\|H\bar{x}^{\frac{a}{2}}\|_{L^2}^2,
\end{equation}
and
\begin{align}\label{10.5}
|\bar{I}_2| & \leq \int|\nabla u||H|^2\bar{x}^adx \notag \\
&\leq \|\nabla u\|_{L^2}\|H\bar{x}^{\frac{a}{2}}\|_{L^4}^2
 \notag \\
&\leq C\|H\bar{x}^{\frac{a}{2}}\|_{L^2}
\|\nabla(H\bar{x}^{\frac{a}{2}})\|_{L^2}
 \notag \\
&\leq C\|H\bar{x}^{\frac{a}{2}}\|_{L^2}
\big(\|\nabla H\bar{x}^{\frac{a}{2}}\|_{L^2}
+\|H\nabla\bar{x}^{\frac{a}{2}}\|_{L^2}\big) \notag \\
&\leq C\|H\bar{x}^{\frac{a}{2}}\|_{L^2}
\big(\|\nabla H\bar{x}^{\frac{a}{2}}\|_{L^2}
+\|H\bar{x}^{\frac{a}{2}}\|_{L^2}
\|\bar{x}^{-1}\nabla\bar{x}\|_{L^\infty}\big)
 \notag \\
&\leq C\|H\bar{x}^{\frac{a}{2}}\|_{L^2}^2
+\frac{\nu}{4}\|\nabla H\bar{x}^{\frac{a}{2}}\|_{L^2}^2,
\end{align}
due to \eqref{2.4}, \eqref{3.5}, and \eqref{z10.6}. Similarly,
it follows from H{\"o}lder's inequality, \eqref{z10.6}, \eqref{2.4}, and \eqref{3.a2} that, for $a>1$,
\begin{align}\label{10.6}
|\bar{I}_3| & \leq C\int|H|^2\bar{x}^a
\bar{x}^{-\frac12}|u|\bar{x}^{-\frac12+\frac15}dx \notag \\
&\leq C\|H\bar{x}^{\frac{a}{2}}\|_{L^4}
\|H\bar{x}^{\frac{a}{2}}\|_{L^2}
\|u\bar{x}^{-\frac34}\|_{L^4}\|\bar{x}^{-\frac{1}{20}}\|_{L^\infty}
\notag \\
&\leq C\|H\bar{x}^{\frac{a}{2}}\|_{L^2}^2
+\frac{\nu}{4}\|\nabla H\bar{x}^{\frac{a}{2}}\|_{L^2}^2.
\end{align}
Putting \eqref{10.4}--\eqref{10.6} into \eqref{lv4.1}, we thus deduce from Gronwall's inequality and \eqref{1.9} that
\begin{align}\label{lbqnew-gj10}
\sup_{0\le t\le T} \|H\bar{x}^{\frac{a}{2}}\|_{L^2}^2
+\int_{0}^{T}\|\nabla H\bar{x}^{\frac{a}{2}}\|_{L^2}^2dt\le C.
\end{align}
This along with \eqref{igj1-2} gives the desired \eqref{igj1}.
\hfill $\Box$

\begin{lemma}\label{le-3}
There exists a positive constant $C$ depending on $T$ such that
\begin{align}\label{gj6}
\sup_{0\le t\le T}\big(\|\sqrt{\rho}u_t\|_{L^2}^2+\|H_t\|_{L^2}^2\big) +\int_{0}^{T}\big(\|\nabla u_t\|_{L^2}^2+\|\nabla H_t\|_{L^2}^2\big)dt \le C.
\end{align}
\end{lemma}
\textit{Proof.}
1. Differentiating $\eqref{mhd}_2$ with respect to $t$ gives that
\begin{align}\label{zb1}
\rho u_{tt}+\rho u\cdot \nabla u_t-\mu\Delta u_t
=-\rho_t(u_t+u\cdot\nabla u)-\rho u_t\cdot\nabla u -\nabla P_t
+\left(H\cdot\nabla H\right)_t.
\end{align}
Multiplying \eqref{zb1} by $u_t$ and integrating the resulting equality by parts over $B_R$, we obtain after using $\eqref{mhd}_1$ and $\eqref{mhd}_5$ that
\begin{align}\label{na8}
& \frac{1}{2}\frac{d}{dt} \int \rho|u_t|^2dx+\mu\int|\nabla u_t|^2 dx \notag \\
& \le C\int \rho|u||u_{t}| \left(|\nabla u_t|+|\nabla u|^{2}+|u||\nabla^{2}u|\right)dx +C\int \rho|u|^{2}|\nabla u ||\nabla u_{t}|dx \notag \\
& \quad+C\int \rho|u_t|^{2}|\nabla u|dx +\int H_t\cdot \nabla H\cdot u_tdx+\int H\cdot \nabla H_t\cdot u_t dx\triangleq \sum_{i=1}^5 \hat{I}_i.
\end{align}
It follows from \eqref {local1}, \eqref{3.a2}, \eqref{2.4}, and \eqref{3.5} that
\begin{align}\label{na2}
\hat{I}_1& \le C \|\sqrt{\rho} u\|_{L^{6}}\|\sqrt{\rho}u_{t}\|_{L^{2}}^{\frac12} \|\sqrt{\rho}u_{t}\|_{L^{6}}^{\frac12}\big(\|\nabla u_{t}\|_{L^{2}}
+\|\nabla u\|_{L^4}^2\big) \notag \\
& \quad +C\|\rho^{\frac14}u\|_{L^{12}}^{2}\|\sqrt{\rho}u_{t}\|_{L^{2}}^{\frac12} \|\sqrt{\rho}u_{t}\|_{L^{6}}^{\frac12} \|\nabla^{2} u \|_{L^{2}} \notag \\
& \le C\|\sqrt{\rho}u_{t}\|_{L^{2}}^{\frac12}
\big(\|\sqrt{\rho}u_{t}\|_{L^{2}}+\|\nabla u_{t}\|_{L^2}\big)^{\frac12}
\big(\|\nabla u_{t}\|_{L^{2}}+\|\nabla^2u\|_{L^2}\big) \notag \\
& \le  \frac{\mu}{6}\|\nabla u_{t}\|_{L^2}^{2}+C\|\sqrt{\rho}u_{t}\|_{L^2}^2+C\|\nabla^2u\|_{L^2}^2.
\end{align}
We infer from H\"older's inequality, \eqref{local1}, \eqref{3.a2}, Sobolev's inequality, and \eqref{3.5} that
\begin{align}\label{5.a3}
\hat{I}_2+\hat{I}_3
& \le C\|\sqrt{\rho}u\|_{L^8}^{2}\|\nabla u\|_{L^4}
\|\nabla u_{t}\|_{L^{2}}+C \| \nabla u\|_{L^{2}}
\|\sqrt{\rho}u_{t}\|_{L^{6}}^{\frac32}\|\sqrt{\rho}u_{t}\|_{L^{2}}^{\frac12} \notag \\
& \le C\|\nabla u\|_{H^1}\|\nabla u_{t}\|_{L^2}
+C\big(\|\sqrt{\rho}u_{t}\|_{L^2}+\|\nabla u_{t}\|_{L^2}\big)^{\frac32}\|\sqrt{\rho}u_{t}\|_{L^2}^{\frac12} \notag \\
& \leq \frac{\mu}{6} \| \nabla u_{t}\|_{L^{2}}^{2}+C\|\sqrt{\rho}u_{t}\|_{L^2}^2+C\|\nabla u\|_{H^1}^2.
\end{align}
It follows from integration by parts, H\"older's inequality, \eqref{3.19}, and \eqref{2.4} that
\begin{align}\label{na3}
\hat{I}_4+\hat{I}_5
& =-\int H_t\cdot\nabla u_t\cdot H dx-\int H\cdot\nabla u_t\cdot H_tdx \notag \\
& \leq 2\|\nabla u_{t}\|_{L^2}\|H\|_{L^4}\|H_t\|_{L^4} \notag \\
& \le \frac{\mu}{6}\|\nabla u_{t}\|_{L^2}^2+C\|H_t\|_{L^2}\|\nabla H_t\|_{L^2}
\notag \\
& \le \frac{\mu}{6}\|\nabla u_{t}\|_{L^2}^2
+\varepsilon\|\nabla H_{t}\|_{L^2}^2+C\|H_t\|_{L^2}^2.
\end{align}
Thus, substituting \eqref{na2}--\eqref{na3} into \eqref{na8}, we obtain that
\begin{align}\label{a4.6}
\frac{d}{dt}\|\sqrt{\rho}u_{t}\|_{L^2}^2+\mu\|\nabla u_{t}\|_{L^2}^2
\le C\|H_t\|_{L^2}^2+\varepsilon\|\nabla H_t\|_{L^{2}}^{2}+C\|\sqrt{\rho}u_{t}\|_{L^2}^2+C\|\nabla u\|_{H^1}^2.
\end{align}

2. Differentiating $\eqref{mhd}_4$ with respect to $t$ shows that
\begin{align}\label{lv4.12}
H_{tt}-H_t\cdot\nabla u-H\cdot\nabla u_t+u_t\cdot\nabla H+u\cdot\nabla H_t=\nu\Delta H_t.
\end{align}
Multiplying \eqref{lv4.12} by $H_t$ and integrating the resulting equality  over $B_R$ yield that
\begin{align}\label{lv4.13}
& \frac12\frac{d}{dt}\int|H_t|^2dx+\nu\int|\nabla H_t|^2dx \notag \\
& =\int(H\cdot\nabla)u_t\cdot H_tdx-\int(u_t\cdot\nabla) H\cdot H_tdx
+\int(H_t\cdot\nabla)u\cdot H_tdx-\int(u\cdot\nabla) H_t\cdot H_tdx \notag \\
& \triangleq S_1+S_2+S_3+S_4.
\end{align}
Integration by parts together with \eqref{3.19}, H{\"o}lder's inequality, \eqref{3.v2}, and \eqref{lbqnew-gj10} leads to
\begin{align}\label{lvbo4.14}
S_1+S_2 & = -\int(H\cdot\nabla)H_t\cdot u_tdx+\int(u_t\cdot\nabla) H_t\cdot Hdx \notag \\
& \le 2\|\nabla H_t\|_{L^2} \||u_t||H|\|_{L^2} \notag \\
&\le \frac{\nu}{4}\|\nabla H_t\|_{L^2}^2
+\frac{4}{\nu}\|u_t\bar{x}^{-\frac{a}{4}}\|_{L^8}^2 \|H \bar{x}^{\frac{a}{2}}\|_{L^2}\|H\|_{L^4}
\notag \\
& \le \frac{\nu}{4}\|\nabla H_t\|_{L^2}^2
+C(\nu,a)\|\sqrt{\rho}u_t\|_{L^2}^2+C(\nu,a)\|\nabla u_t\|_{L^2}^2.
\end{align}
By virtue of H{\"o}lder's inequality, \eqref{2.4}, and \eqref{3.5} one has
\begin{align}
S_3 &  \leq \|H_t\|_{L^4}^2\|\nabla u\|_{L^2}
\leq  C \|H_t\|_{L^2}\|\nabla H_t\|_{L^2} \|\nabla u\|_{L^2}
\leq \frac{\nu}{4}\|\nabla H_t\|_{L^2}^2+C\|H_t\|_{L^2}^2.
\end{align}
We derive from integration by parts and $\eqref{mhd}_{5}$ that
\begin{align*}
S_4 =\int(u\cdot\nabla) H_t\cdot H_tdx=-S_4,
\end{align*}
that is
\begin{align}\label{lv4.14}
S_4 =0.
\end{align}
Inserting \eqref{lvbo4.14}--\eqref{lv4.14} into \eqref{lv4.13}, we get
\begin{align}\label{ilv4.14}
& \frac{d}{dt}\|H_t\|_{L^2}^2+\nu\|\nabla H_t\|_{L^2}^2\le  C\|H_t\|_{L^2}^2+C\|\sqrt{\rho}u_t\|_{L^2}^2
+C_2\|\nabla u_t\|_{L^2}^2
\end{align}
for some positive constant $C_2$ depending on $\nu$ and $a$. Adding \eqref{ilv4.14} multiplied by $\frac{\mu}{2C_2}$ to \eqref{a4.6} and then choosing $\varepsilon=\frac{\mu\nu}{4C_2}$, we arrive at
\begin{align}\label{3.50}
& \frac{d}{dt}\Big(\|\sqrt{\rho}u_{t}\|_{L^2}^2+\frac{\mu}{2C_2}\|H_t\|_{L^2}^2\Big)
+\frac{\mu}{2}\|\nabla u_{t}\|_{L^2}^2+\frac{\mu\nu}{4C_2}\|\nabla H_t\|_{L^2}^2 \notag \\
& \le C\big(\|\sqrt{\rho}u_{t}\|_{L^2}^2+\|H_t\|_{L^2}^2\big)
+C\|\nabla u\|_{H^1}^2.
\end{align}

3. It follows from \eqref{2.5}, Young's inequality, \eqref{z10.6}, \eqref{3.v2}, \eqref{2.4}, \eqref{3.4}, and \eqref{3.5} that
\begin{align}\label{3.22}
\|u\bar x^{-\frac{a}{2}}\|_{L^\infty}
& \le C\|\nabla(u\bar x^{-\frac{a}{2}})\|_{L^3}^{\frac35}
\|u\bar x^{-\frac{a}{2}}\|_{L^4}^{\frac25} \notag \\
& \le C\big(\|u\bar x^{-\frac{a}{2}}\|_{L^4}
+\|\nabla(u\bar x^{-\frac{a}{2}})\|_{L^3}\big) \notag \\
& \le C\big(\|u\bar x^{-\frac{a}{2}}\|_{L^4}
+\|u\bar x^{-\frac{a}{2}-1+\frac{a}{2}}\|_{L^3}
+\|\nabla u\|_{L^3}\big)
\notag \\
& \le C\Big(\|\sqrt{\rho}u\|_{L^2}+\|\nabla u\|_{L^2}
+\|\nabla u\|_{L^2}^{\frac23}\|\nabla u\|_{H^1}^{\frac13}\Big) \notag \\
& \leq C\|\nabla^2u\|_{L^2}+C.
\end{align}
We deduce from \eqref{mhd}$_4$, \eqref{3.19}, and \eqref{2.4} that
\begin{align*}
\|H_t\|_{L^2}^2& \leq C\big(\|\Delta H\|_{L^2}^2+\|u\cdot \nabla H\|_{L^2}^2+\|H\cdot \nabla u\|_{L^2}^2\big) \notag \\
& \le C\big(\|\nabla^2 H\|_{L^2}^2+\|u \bar x^{-\frac{a}{2}}\|_{L^{\infty}}^2\|\nabla H \bar x^{\frac{a}{2}}\|_{L^2}^2
+\|H\|_{L^4}^2\|\nabla u \|_{L^4}^2\big) \notag \\
& \le C\big(\|\nabla^2 H\|_{L^2}^2+\|u \bar x^{-\frac{a}{2}}\|_{L^{\infty}}^2\|\nabla H \bar x^{\frac{a}{2}}\|_{L^2}^2
+\|\nabla u \|_{L^2}\|\nabla u \|_{H^1}\big),
\end{align*}
which together with \eqref{3.22} and \eqref{1.9} yields that
\begin{align}\label{htgj2}
\|H_t(0)\|_{L^2}^2\le C.
\end{align}
From \eqref{2.5} and \eqref{3.22}, one has
\begin{align}\label{3.53}
\|\sqrt{\rho} u\|_{L^\infty}
& \le\|\rho \bar{x}^a\|_{L^{\infty}}^\frac12\|u\bar{x}^{-\frac{a}{2}}\|_{L^\infty} \notag \\
& \le C\|\rho \bar{x}^a\|^{\frac{q-2}{4(q-1)}}_{L^2}
\|\nabla(\rho \bar{x}^a)\|^{\frac{q}{4(q-1)}}_{L^q}
\big(\|\nabla^2u\|_{L^2}+1\big),
\end{align}
which combined with \eqref{mhd}$_2$, \eqref{c.c}, and \eqref{1.9} leads to
\begin{align*}
\int \rho|u_t|^2(x,0)dx
& \le\lim_{t \rightarrow 0}\sup\int\rho^{-1}|\mu\Delta u+H\cdot\nabla H
-\nabla P-\rho u\cdot\nabla u|^2dx \\
& \le 2\|g_1\|_{L^2}^2+2\|\sqrt{\rho}u(0)\|_{L^\infty}^2\|\nabla u(0)\|_{L^2}^2 \le C.
\end{align*}
This along with \eqref{3.50}, Gronwall's inequality, \eqref{htgj2}, \eqref{3.4}, and \eqref{3.5} leads to \eqref{gj6}. \hfill $\Box$

\begin{lemma} \label{le-3'}
Let $q$ be as in Theorem \ref{t1}, then there exists a positive constant $C$ depending on $T$ such that
\begin{align}\label{gj10'}
&\sup_{0\le t\le T}\big(\|\nabla^2 u\|_{L^2}^2+\|\nabla P\|_{L^2}^2
+\|\nabla^2 H\|_{L^2}^2+\|\nabla H\bar{x}^{\frac{a}{2}}\|_{L^2}^2\big) \notag \\
& \quad +\int_{0}^{T}\Big(\|\nabla^2u\|_{ L^q}^{\frac{q+1}{q}}
+\|\nabla P\|_{L^q}^{\frac{q+1}{q}}+\|\nabla^2u\|_{L^q}^2
+\|\nabla P\|_{ L^q}^2+\|\nabla^2H\bar{x}^{\frac{a}{2}}\|_{L^2}^2\Big)dt\le C.
\end{align}
\end{lemma}
\textit{Proof.}
1. Multiplying \eqref{mhd}$_4$ by $\Delta H\bar{x}^a$ and integrating by parts lead to
\begin{align}\label{AMSS5}
&\frac{1}{2}\frac{d}{dt}\int |\nabla H|^2\bar{x}^adx+\nu \int |\Delta H|^2\bar{x}^adx \notag \\
& \le C\int|\nabla H| |H| |\nabla u| |\nabla\bar{x}^a|dx+C\int|\nabla H|^2|u| |\nabla\bar{x}^a|dx+C\int|\nabla H| |\Delta H| |\nabla \bar{x}^a|dx \notag \\
& \quad+ C\int |H||\nabla u||\Delta H|\bar{x}^adx+C\int |\nabla u||\nabla H|^2 \bar{x}^adx
\triangleq \sum_{i=1}^5 J_i.
\end{align}
Using \eqref{lbqnew-gj10}, \eqref{3.a2}, H\"older's inequality, \eqref{2.4}, and \eqref{3.5}, we get by some direct calculations that
\begin{align*}
J_1\le & C\|H\bar{x}^{\frac{a}{2}}\|_{L^4}\|\nabla u\|_{L^4}\|\nabla H\bar{x}^{\frac{a}{2}}\|_{L^2}\\
\le & C\|H\bar{x}^{\frac{a}{2}}\|_{L^2}^{\frac12}
\big(\|\nabla H\bar{x}^{\frac{a}{2}}\|_{L^2}
+\|H\bar{x}^{\frac{a}{2}}\|_{L^2}\big)^{\frac12}
\|\nabla u\|_{L^2}^{\frac12}\|\nabla^2u\|_{L^2}^{\frac12}
\|\nabla H\bar{x}^{\frac{a}{2}}\|_{L^2}\\
\le & C\|\nabla H\bar{x}^{\frac{a}{2}}\|_{L^2}^2+C\|\nabla u\|_{H^1}^2,\\
J_2\leq & C\||\nabla H|^{2-\frac{2}{3a}}\bar{x}^{a-\frac13}\|_{L^{\frac{6a}{6a-2}}} \|u\bar{x}^{-\frac13}\|_{L^{6a}}\||\nabla H|^{\frac{2}{3a} }\|_{L^{6a}}\\
\le & C\|\nabla H\bar{x}^{\frac{a}{2}} \|_{L^2}^\frac{6a-2}{3a}\|\nabla H\|_{L^4}^\frac{2}{3a} \\
\leq & C\|\nabla H\bar{x}^{\frac{a}{2}}\|_{L^2}^2
+C\|\nabla H\|_{L^4}^2\\
\leq & C\|\nabla H \bar{x}^{\frac{a}{2}}\|_{L^2}^2
+\frac{\nu}{4}\|\Delta H \bar{x}^{\frac{a}{2}}\|_{L^2}^2,\\
J_3+J_4\le &\frac{\nu}{4}\|\Delta H\bar{x}^{\frac{a}{2}}\|_{L^2}^2
+C\|\nabla H\bar{x}^{\frac{a}{2}}\|_{L^2}^2
+C\|H\bar{x}^{\frac{a}{2}}\|_{L^4}^2 \|\nabla u\|_{L^4}^2\\
\le &\frac{\nu}{4}\|\Delta H\bar{x}^{\frac{a}{2}}\|_{L^2}^2+C\|\nabla H\bar{x}^{\frac{a}{2}}\|_{L^2}^2
+C\|H\bar{x}^{\frac{a}{2}}\|_{L^2}\big(\|\nabla H\bar{x}^{\frac{a}{2}}\|_{L^2}+\|H\bar{x}^{\frac{a}{2}}\|_{L^2}\big)
\|\nabla u\|_{L^2} \|\nabla u\|_{H^1} \\
\le & \frac{\nu}{4}\|\Delta H\bar{x}^{\frac{a}{2}}\|_{L^2}^2
+C\|\nabla H\bar{x}^{\frac{a}{2}}\|_{L^2}^2+C\|\nabla u\|_{H^1}^2,\\
J_5 \le & C\|\nabla u\|_{L^\infty} \|\nabla H\bar{x}^{\frac{a}{2}}\|_{L^2}^2
\le  C\Big(1+\|\nabla^2u\|_{L^q}^{\frac{q+1}{q}}\Big)\|\nabla H\bar{x}^{\frac{a}{2}}\|_{L^2}^2.
\end{align*}
Substituting the above estimates into \eqref{AMSS5}
and noting the following fact
\begin{align*}
\int|\nabla^2H|^2\bar{x}^adx
& = \int|\Delta H|^2\bar{x}^adx
-\int\partial_i\partial_kH\cdot\partial_kH\partial_i\bar{x}^adx
+\int\partial_i\partial_iH\cdot\partial_kH\partial_k\bar{x}^adx \\
& \leq \int|\Delta H|^2\bar{x}^adx+\frac12\int|\nabla^2H|^2\bar{x}^adx
+C\int|\nabla H|^2\bar{x}^adx,
\end{align*}
we derive that
\begin{align}\label{AMSS10}
& \frac{d}{dt}\|\nabla H\bar{x}^{\frac{a}{2}}\|_{L^2}^2
+\nu\|\nabla^2H\bar{x}^{\frac{a}{2}}\|_{L^2}^2
\le C\Big(1+\|\nabla^2 u\|_{L^q}^{\frac{q+1}{q}}\Big)
\|\nabla H\bar{x}^{\frac{a}{2}}\|_{L^2}^2+C\|\nabla u\|_{H^1}^2.
\end{align}
Now we claim that
\begin{align}\label{gj7}
\int_{0}^{T}\Big(\|\nabla^2u\|_{ L^q}^{\frac{q+1}{q}}
+\|\nabla P\|_{L^q}^{\frac{q+1}{q}}+\|\nabla^2u\|_{L^q}^2
+\|\nabla P\|_{ L^q}^2\Big)dt\le C,
\end{align}
whose proof will be given at the end of this proof. Thus, we infer from \eqref{AMSS10}, Gronwall's inequality, \eqref{gj7}, \eqref{3.4}, and \eqref{3.5} that
\begin{align}\label{igj10'}
\sup_{0\le t\le T}\|\nabla H\bar{x}^{\frac{a}{2}}\|_{L^2}^2
+\int_{0}^{T}\|\nabla^2H\bar{x}^{\frac{a}{2}}\|_{L^2}^2dt
\le C.
\end{align}

2. It deduces from $\eqref{mhd}_4$, the standard $L^2$-estimate of elliptic equations, \eqref{3.a2}, H\"older's inequality, \eqref{3.5}, \eqref{3.19}, and Gagliardo-Nirenberg inequality that
\begin{align}\label{AMSS11}
\|\nabla^2H\|^2_{L^2}
&\leq C\|H_t\|^2_{L^2}+C\||u||\nabla H|\|_{L^2}^2+C\||H||\nabla u|\|^2_{L^2} \notag \\
& \leq C\|H_t\|^2_{L^2}+C\|u \bar{x}^{-\frac{a}{4}}\|_{L^8}^2\|\nabla H \bar{x}^{\frac{a}{2}}\|_{L^2}\|\nabla H\|_{L^4}
+C\|H\|_{L^4}^2\|\nabla u\|_{L^4}^2 \notag \\
& \leq C\|H_t\|^2_{L^2}+C\|\nabla H \bar{x}^{\frac{a}{2}}\|_{L^2}^2
+C\|\nabla H\|_{L^2}\|\nabla H\|_{H^1}
+C\|\nabla u\|_{L^2}\|\nabla u\|_{H^1} \notag \\
& \leq C\|H_t\|^2_{L^2}+C\|\nabla H \bar{x}^{\frac{a}{2}}\|_{L^2}^2
+\frac14\|\nabla^2 H\|_{L^2}^2+\frac14\|\nabla^2 u\|_{L^2}^2+C.
\end{align}
It follows from \eqref{3.2}, \eqref{3.a2}, \eqref{3.19}, \eqref{3.5}, and \eqref{2.4} that
\begin{align*}
\|\nabla^2u\|_{L^2}^2+\|\nabla P\|_{L^2}^2
&\leq C\|\rho u_t\|_{L^2}^2+C\|\rho u\cdot\nabla u\|_{L^2}^2
+C\||H||\nabla H|\|_{L^2}^2
\notag \\
&\leq C \|\sqrt{\rho}u_t\|_{L^2}^2
+C\|\rho u\|_{L^4}^2\|\nabla u\|_{L^4}^2+C\|H\|_{L^4}^2\|\nabla H\|_{L^4}^2 \notag \\
&\leq C \|\sqrt{\rho}u_t\|_{L^2}^2
+C\|\nabla u\|_{L^2}\|\nabla u\|_{H^1}+C\|\nabla H\|_{L^2}\|\nabla H\|_{H^1} \notag \\
&\leq C \|\sqrt{\rho}u_t\|_{L^2}^2
+\frac14\|\nabla^2u\|_{L^2}^2+\frac14\|\nabla^2H\|_{L^2}^2+C,
\end{align*}
which combined with \eqref{AMSS11}, \eqref{gj6}, and \eqref{igj10'} yields that
\begin{align}\label{iAMSS12}
\sup_{0\le t\le T}\big(\|\nabla^2u\|_{L^2}^2+\|\nabla P\|_{L^2}^2
+\|\nabla^2 H\|^2_{L^2}\big)\leq C.
\end{align}

3. To finish the proof of Lemma \ref{le-3'}, it suffices to show \eqref{gj7}. Choosing $p=q$ in \eqref{3.16}, we derive from \eqref{3.1}, H{\"o}lder's inequality, \eqref{3.a2}, \eqref{2.4}, and \eqref{3.5} that
\begin{align}\label{lv3.60}
\|\nabla^2u\|_{L^q}+\|\nabla P\|_{L^q}
& \le C\big(\|\rho u_t\|_{L^q}+\|\rho u\cdot\nabla u\|_{L^q}+\||H||\nabla H|\|_{L^q} \big)\notag \\
&\le C\|\sqrt{\rho}u_t\|_{L^2}^{\frac{2(q-1)}{q^2-2}}
\|\sqrt{\rho}u_t\|_{L^{q^2}}^{\frac{q^2-2q}{q^2-2}}
+C\|\rho u\|_{L^{2q}}\|\nabla u\|_{L^{2q}}
+C\|H\|_{L^{2q}}\|\nabla H\|_{L^{2q}} \notag \\
& \le C\|\sqrt{\rho}u_t\|_{L^2}^{\frac{2(q-1)}{q^2-2}}
\big(\|\sqrt{\rho}u_t\|_{L^2}+\|\nabla u_t\|_{L^2}\big)^{\frac{q^2-2q}{q^2-2}}
+\|\nabla u\|_{L^2}^{\frac1q}\|\nabla^2u\|_{L^2}^{\frac{q-1}{q}} \notag \\
& \quad +C\|H\|_{L^2}^{\frac1q}\|\nabla H\|_{L^2}^{\frac{q-1}{q}}
\|\nabla H\|_{L^2}^{\frac1q}\|\nabla^2H\|_{L^2}^{\frac{q-1}{q}} \notag \\
&\le C\|\nabla u_t\|_{L^{2}}^{\frac{q^2-2q}{q^2-2}}
+C\|\nabla^2u\|_{L^2}^{\frac{q-1}{q}}
+C\|\nabla^2H\|_{L^2}^{\frac{q-1}{q}},
\end{align}
which together with Young's inequality, \eqref{3.5}, and \eqref{gj6} implies that
\begin{align} \label{olv3.61}
\int_0^T\Big(\|\nabla^2u\|_{L^q}^{\frac{q+1}{q}}
+\|\nabla P\|_{L^q}^{\frac{q+1}{q}}\Big)dt
& \le C\int_0^T
\left(\|\nabla u_t\|_{L^2}^{\frac{q^2-q-2}{q^2-2}}
+\|\nabla^2u\|_{L^2}^{\frac{q^2-1}{q^2}}
+\|\nabla^2H\|_{L^2}^{\frac{q^2-1}{q^2}}\right)dt
 \notag \\
& \le C\int_0^T\left(\|\nabla u_t\|_{L^2}^2
+\|\nabla^2u\|_{L^2}^2+\|\nabla^2H\|_{L^2}^2+1\right)dt\le C(T),
\end{align}
and
\begin{align}\label{olv3.62}
\int_0^T\big(\|\nabla^2u\|_{L^q}^2+\|\nabla P\|_{L^q}^2\big)dt
& \le C\int_0^T
\left[\big(\|\nabla u_t\|_{L^{2}}^2\big)^{\frac{q^2-2q}{q^2-2}}
+(\|\nabla^2u\|_{L^2}^2)^{\frac{q-1}{q}}
+(\|\nabla^2H\|_{L^2}^2)^{\frac{q-1}{q}}\right]dt
 \notag \\
& \le C\int_0^T\left(\|\nabla u_t\|_{L^2}^2
+\|\nabla^2u\|_{L^2}^2+\|\nabla^2H\|_{L^2}^2+1\right)dt\le C(T),
\end{align}
where we have used $\frac{q^2-2q}{q^2-2}\in(0,1)$ due to $q>2$.
One thus obtains \eqref{gj7} from \eqref{olv3.61} and \eqref{olv3.62}.
\hfill $\Box$

\begin{lemma}\label{le-4}
Let $q$ be as in Theorem \ref{t1}, then there exists a positive constant $C$ depending on $T$ such that
\begin{align}\label{gj8}
\sup\limits_{0\le t\le T}
\big(\|\rho\bar x^a\|_{H^1\cap W^{1,q}}+\|\rho_t\|_{L^2\cap L^q}\big)\le C.
\end{align}
\end{lemma}
\textit{Proof.}
1. We derive from \eqref{mhd}$_1$ and \eqref{mhd}$_5$ that $\rho\bar{x}^a$ satisfies
\begin{equation}\label{A}
\partial_{t}(\rho\bar{x}^a)+ u \cdot\nabla(\rho\bar{x}^a)
-a\rho\bar{x}^{a}u\cdot\nabla(\ln\bar{x})=0.
\end{equation}
Multiplying \eqref{A} by $(\rho\bar{x}^a)^{r-1}$ for $r\in[2,q]$ and then integrating the resultant equation over $\mathbb{R}^2$, we deduce that
\begin{align*}
\frac{d}{dt}\int(\rho\bar{x}^a)^rdx
=ar\int (\rho\bar{x}^a)^ru\cdot\nabla(\ln\bar{x})dx,
\end{align*}
which leads to
\begin{align}\label{7.2}
\frac{d}{dt}\|\rho\bar{x}^a\|_{L^r}^r \leq
& ar\|u\cdot\nabla(\ln\bar{x})\|_{L^\infty}\|\rho\bar{x}^a\|_{L^r}^r
\end{align}
Similarly to \eqref{3.22}, we obtain after using \eqref{z10.6} and \eqref{iAMSS12} that
\begin{align}\label{7.3}
\|u\cdot\nabla(\ln\bar{x})\|_{L^\infty}
=\|u\cdot\bar{x}^{-1}\nabla\bar{x}\|_{L^\infty}
\leq C.
\end{align}
This along with \eqref{7.2} and Gronwall's inequality yields that
\begin{align}\label{7.1}
\sup\limits_{0\le t\le T} \|\rho\bar x^a\|_{L^r}\le C.
\end{align}

2. Operating $\nabla$ to \eqref{A} and then multiplying the resultant equation by $|\nabla(\rho\bar{x}^a)|^{r-2}\nabla(\rho\bar{x}^a)$ for $r\in[2,q]$,
we obtain after integration by parts that
\begin{align}\label{6.4}
& \frac{d}{dt}\|\nabla(\rho\bar{x}^a)\|_{L^r} \notag \\
& \leq C\big(1+\|\nabla u \|_{L^\infty}
+\|u\cdot\nabla(\ln\bar{x})\|_{L^\infty}\big)
\|\nabla(\rho\bar{x}^a)\|_{L^r}
+C\|\rho\bar{x}^a\|_{L^\infty}\big(\||\nabla u ||\nabla(\ln\bar{x})|\|_{L^r}+\||u||\nabla^{2}(\ln\bar{x})|\|_{L^r}\big) \notag \\
& \leq C\big(1+\|\nabla u \|_{L^\infty}\big)
\|\nabla(\rho\bar{x}^a)\|_{L^r}
+C\|\rho\bar{x}^a\|_{L^\infty}\big(\||\nabla u ||\nabla(\ln\bar{x})|\|_{L^r}+\||u||\nabla^{2}(\ln\bar{x})|\|_{L^r}\big),
\end{align}
due to \eqref{7.3}.
By \eqref{2.5}, \eqref{3.5}, and Young's inequality, we see that
\begin{align*}
\|\nabla u\|_{L^\infty}\leq C\|\nabla u\|_{L^2}^{\frac{q-2}{2q-2}}\|\nabla^2u\|_{L^q}^{\frac{q}{2q-2}}
\leq C\|\nabla^2u\|_{L^q}^2+C.
\end{align*}
From \eqref{2.5}, \eqref{7.1}, and Young's inequality, we have
\begin{align*}
\|\rho\bar{x}^a\|_{L^\infty}
\leq C\|\rho\bar{x}^a\|_{L^2}^{\frac{q-2}{2q-2}}
\|\nabla(\rho\bar{x}^a)\|_{L^q}^{\frac{q}{2q-2}}
\leq C\|\nabla(\rho\bar{x}^a)\|_{L^q}^{\frac{q}{2q-2}}
\leq C+C\|\nabla(\rho\bar{x}^a)\|_{L^q}.
\end{align*}
Applying \eqref{z10.6} and \eqref{2.4}, we get from \eqref{3.5} and \eqref{iAMSS12} that
\begin{align*}
\||\nabla u ||\nabla(\ln\bar{x})|\|_{L^r}
\leq C\|\nabla u\|_{L^r}\|\bar{x}^{-\frac{4+a}{8+a}}\|_{L^\infty}
\leq C\|\nabla u\|_{L^2}^{\frac2r}\|\nabla u\|_{H^1}^{\frac{r-2}{r}}
\leq C.
\end{align*}
Moreover, it follows from \eqref{z10.6} and \eqref{3.a2} that
\begin{align*}
\||u||\nabla^{2}(\ln\bar{x})|\|_{L^r}\leq C.
\end{align*}
As a consequence, inserting the above estimates into \eqref{6.4}, we derive that
\begin{align}\label{7.4}
\frac{d}{dt}\|\nabla(\rho\bar{x}^a)\|_{L^r}
\leq C\big(1+\|\nabla^2 u\|_{L^q}^2\big)
\big(1+\|\nabla(\rho\bar{x}^a)\|_{L^r}+\|\nabla(\rho\bar{x}^a)\|_{L^q}\big).
\end{align}
Hence, choosing $r=q$, then we obtain from Gronwall's inequality and \eqref{olv3.62} that
\begin{align}\label{3.68}
\sup\limits_{0\le t\le T} \|\nabla(\rho\bar x^a)\|_{L^q}\le C.
\end{align}
This along with \eqref{7.4}, Gronwall's inequality, and \eqref{olv3.62} implies that
\begin{align*}
\sup\limits_{0\le t\le T}\|\nabla(\rho\bar x^a)\|_{L^2}\le C,
\end{align*}
which together with \eqref{7.1} and \eqref{3.68} gives that
\begin{align}\label{x11}
\sup\limits_{0\le t\le T}\|\rho\bar x^a\|_{H^1\cap W^{1,q}}\le C.
\end{align}

3. It follows from \eqref{3.22} and \eqref{iAMSS12} that
\begin{align}\label{x12}
\sup\limits_{0\le t\le T}\|u\bar x^{-\frac{a}{2}}\|_{L^\infty}\le C.
\end{align}
Noticing that
\begin{align*}
\rho_t=-u\cdot\nabla\rho
=-u\bar x^{-\frac{a}{2}}\cdot\nabla\rho\bar x^a
\bar x^{-\frac{a}{2}},
\end{align*}
which combined with \eqref{x12} and \eqref{x11} yields that
\begin{align}\label{x13}
\sup\limits_{0\le t\le T}\|\rho_t\|_{L^2\cap L^q}\le C.
\end{align}
This completes the proof of Lemma \ref{le-4}.
\hfill $\Box$

\begin{lemma}\label{lemma 3.8}
Let $q$ be as in Theorem \ref{t1}, then there exists a positive constant $C$ depending on $T$ such that
\begin{align}\label{wdgj}
&\sup\limits_{0\le t\le T} \left(\|\sqrt{\rho} \theta\|_{L^2}^2+\|\nabla \theta\|_{H^1}^2+\|\sqrt{\rho} \theta_t\|_{L^2}^2\right) \notag \\
& \quad +\int_0^T\Big(\|\sqrt{\rho} \theta_t\|_{L^2}^2+\|\nabla^2 \theta\|_{L^2}^2+\|\nabla^2 \theta\|_{L^q}^{\frac{q+1}{q}}+\|\nabla^2\theta\|_{L^q}^2
+\|\nabla\theta_t\|_{L^2}^2\Big)dt\le  C.
\end{align}
\end{lemma}
\textit{Proof.}
1. Choosing $b_1\le\frac{a}{2}$ in Lemma \ref{lemma2.5}, then for $0<b<\min\{b_1,1\}$, we have
$$\bar{x}^b\le C\big(1+|x|^{b_1}\big)<C\bar{x}^{\frac{a}{2}}.$$
Thus it follows from Lemma \ref{lemma2.5}, \eqref{3.4}, \eqref{3.5}, and \eqref{gj8} that
\begin{align}\label{3.8.1}
& \int \Big[\frac{\mu}{2}|\nabla u+(\nabla u)^{tr}|^{2}+\nu(\curl H)^{2}\Big]\bar{x}^b dx \notag \\
&\le C\|\nabla u\|_{L^2}^2+C\|\nabla H\|_{L^2}^2+C\int\big(\rho|\theta_t|+\rho|u||\nabla \theta|\big)|x|^{b_1}dx    \notag \\
&\le C+C\|\sqrt{\rho} \bar{x}^{\frac{a}{2}}\|_{L^2\bigcap L^\infty}\big(\|\sqrt{\rho}\theta_t\|_{L^2}
+\|\sqrt{\rho} u\|_{L^2}\|\nabla\theta\|_{L^2}\big) \notag \\
& \le C+C\big(\|\sqrt{\rho}\theta_t\|_{L^2}+\|\nabla\theta\|_{L^2}\big).
\end{align}
Multiplying \eqref{mhd}$_3$ by $\theta$ and integration by parts, one has that
\begin{align}\label{3.8.2}
\frac{c_v}{2}\int\rho\theta^2dx+\kappa \int|\nabla \theta|^2dx
= \int \Big[\frac{\mu}{2}|\nabla u+(\nabla u)^{tr}|^{2}+\nu(\curl H)^{2}\Big]\theta dx.
\end{align}
For simplicity, setting $Z\triangleq\Big[\frac{\mu}{2}|\nabla u+(\nabla u)^{tr}|^{2}+\nu(\curl H)^{2}\Big]$, then we infer from \eqref{3.v2}, \eqref{3.8.1}, \eqref{3.4}, and \eqref{iAMSS12} that
\begin{align}\label{3.8.3}
& \int \Big[\frac{\mu}{2}|\nabla u+(\nabla u)^{tr}|^{2}+\nu(\curl H)^{2}\Big]\theta dx \notag \\
& \le C\|\theta \bar{x}^{-\frac{b}{2}}\|_{L^6}
\|\sqrt{Z}\bar{x}^{\frac{b}{2}}\|_{L^2}
\big(\|\nabla u\|_{L^3}+\|\nabla H\|_{L^3}\big) \notag \\
& \le C\big(\|\sqrt{\rho}\theta\|_{L^2}+\|\nabla\theta\|_{L^2}\big)
\big(\|\sqrt{\rho}\theta_t\|_{L^2}+\|\nabla\theta\|_{L^2}\big)^{\frac12} \notag \\
& \le \frac{\kappa}{2}\|\nabla\theta\|_{L^2}^2+C\|\sqrt{\rho}\theta\|_{L^2}^2
+C\|\sqrt{\rho}\theta_t\|_{L^2}^2+C,
\end{align}
due to $\sqrt{Z}\leq C(|\nabla u|+|\nabla H|)$.
Inserting \eqref{3.8.3} into \eqref{3.8.2}, one obtains that
\begin{align}\label{3.8.5}
c_v\frac{d}{dt}\|\sqrt{\rho} \theta\|_{L^2}^2
+\kappa\|\nabla\theta\|_{L^2}^2
\le C\|\sqrt{\rho}\theta\|_{L^2}^2
+C\|\sqrt{\rho}\theta_t\|_{L^2}^2+C.
\end{align}

2. Multiplying \eqref{mhd}$_3$ by $\theta_t$ gives that
\begin{align}\label{3.8.6}
\frac{\kappa}{2}\frac{d}{dt}\|\nabla\theta\|_{L^2}^2
+c_v\|\sqrt{\rho}\theta_t\|_{L^2}^2=-c_v\int \rho u\cdot\nabla\theta \theta_tdx
+\int Z\theta_tdx.
\end{align}
By virtue of H{\"o}lder's inequality, \eqref{3.a2}, \eqref{3.v2}, and \eqref{7.1}, one has that
\begin{align}\label{3.76}
-c_v\int \rho u \cdot\nabla\theta\theta_t dx
& \le c_v\|\rho \bar{x}^a\|_{L^q}\|u \bar{x}^{-\frac{a}{2}}\|_{L^{\frac{4q}{q-2}}}\|\theta_t \bar{x}^{-\frac{a}{2}}\|_{L^{\frac{4q}{q-2}}}\|\nabla\theta\|_{L^2} \notag \\
& \le C\big(\|\sqrt{\rho}\theta_t\|_{L^2}
+\|\nabla\theta_t\|_{L^2}\big)\|\nabla\theta\|_{L^2} \notag \\
& \le \frac{\kappa}{8}\|\nabla\theta_t\|_{L^2}^2
+\frac{c_v}{4}\|\sqrt{\rho}\theta_t\|_{L^2}^2+C\|\nabla\theta\|_{L^2}^2.
\end{align}
We deduce from H{\"o}lder's inequality, \eqref{3.v2}, \eqref{3.8.1}, \eqref{3.4}, and \eqref{iAMSS12} that
\begin{align}\label{3.77}
\int Z\theta_t dx
& \le C\|\theta_t \bar{x}^{-\frac{b}{2}}\|_{L^6}
\|\sqrt{Z}\bar{x}^\frac{b}{2}\|_{L^2}\big(\|\nabla u\|_{L^3}+\|\nabla H\|_{L^3}\big)
\notag \\
& \le C\big(\|\sqrt{\rho}\theta_t\|_{L^2}+\|\nabla\theta_t\|_{L^2}\big)
\big(\|\sqrt{\rho}\theta_t\|_{L^2}+\|\nabla\theta\|_{L^2}\big)^{\frac12} \notag \\
& \le \frac{\kappa}{8}\|\nabla\theta_t\|_{L^2}^2
+\frac{c_v}{4}\|\sqrt{\rho}\theta_t\|_{L^2}^2+C\|\nabla\theta\|_{L^2}^2+C.
\end{align}
Substituting \eqref{3.76} and \eqref{3.77} into \eqref{3.8.6} leads to
\begin{align}\label{3.8.7}
\frac{\kappa}{2}\frac{d}{dt}\|\nabla\theta\|_{L^2}^2
+\frac{c_v}{2}\|\sqrt{\rho}\theta_t\|_{L^2}^2
\leq\frac{\kappa}{4}\|\nabla \theta_t\|_{L^2}^2
+C\|\nabla\theta\|_{L^2}^2+C.
\end{align}

3. Differentiating \eqref{mhd}$_3$ with respect to $t$ and multiplying the resulting equation by $\theta_t$ yield that
\begin{align}\label{3.8.8}
\frac{c_v}{2}\frac{d}{dt}\|\sqrt{\rho}\theta_t\|_{L^2}^2+\kappa\|\nabla \theta_t\|_{L^2}^2
=-c_v\int\rho_t|\theta_t|^2 dx-c_v\int(\rho u)_t\cdot\nabla \theta\theta_t dx
+\int Z_t\theta_tdx \triangleq \sum_{i=1}^{3}L_i.
\end{align}
It follows from \eqref{mhd}$_1$ and integration by parts that
\begin{align*}
L_1&=-c_v\int\rho_t|\theta_t|^2dx\notag\\
&=-2c_v\int \rho u\cdot \nabla\theta_t \theta_tdx\notag\\
&\le\frac{\kappa}{12}\|\nabla \theta_t\|_{L^2}^2+C\|\sqrt{\rho} u\|_{L^\infty}^2\|\sqrt{\rho}\theta_t\|_{L^2}^2\notag\\
&\le\frac{\kappa}{12}\|\nabla \theta_t\|_{L^2}^2+C\|\sqrt{\rho}\theta_t\|_{L^2}^2,
\end{align*}
where we have used
\begin{align}\label{ZZ}
\|\sqrt{\rho} u\|_{L^\infty}\leq C,
\end{align}
due to \eqref{3.53} and \eqref{iAMSS12}.
In view of \eqref{mhd}$_1$ and \eqref{ZZ}, we obtain from integration by parts that
\begin{align*}
L_2 & =-\int(\rho u)_t\nabla \theta\theta_t dx\notag\\
&=-\int\rho u\cdot\nabla(\theta_t\nabla\theta)dx \notag\\
& \le \|\sqrt{\rho}u\|_{L^\infty}
\|\sqrt{\rho}\theta_t\|_{L^2}\|\nabla^2\theta\|_{L^2}
+\|\rho\|_{L^\infty}^{\frac12}\|\sqrt{\rho}u\|_{L^\infty}
\|\nabla\theta_t\|_{L^2}\|\nabla\theta\|_{L^2}
\notag\\
&\le\frac{\kappa}{12}\|\nabla\theta_t\|_{L^2}^2+C\|\sqrt{\rho}\theta_t\|_{L^2}^2
+C\|\nabla\theta\|_{L^2}^2+C\|\nabla^2\theta\|_{L^2}^2.
\end{align*}
Direct calculation gives that
\begin{align*}
Z_t \leq C\sqrt{Z}\big(|\nabla u_t|+|\nabla H_t|\big),
\end{align*}
which combined with H{\"o}lder's inequality,
\eqref{3.v2}, \eqref{3.8.1}, \eqref{3.5}, and \eqref{iAMSS12}
 ensures that
\begin{align}\label{x10}
L_3 & \leq C\int |\theta_t|\sqrt{Z}\big(|\nabla u_t|+|\nabla H_t|\big)dx \notag \\
& \le C\|\theta_t\bar{x}^{-\frac{b}{4}}\|_{L^8}
\|Z^{\frac14}\bar{x}^\frac{b}{4}\|_{L^4}\|Z^{\frac14}\|_{L^8}
\||\nabla u_t|+|\nabla H_t|\|_{L^2}\notag\\
& \le C\|\theta_t\bar{x}^{-\frac{b}{4}}\|_{L^8}
\|\sqrt{Z}\bar{x}^\frac{b}{2}\|_{L^2}^{\frac12}
\big(\|\nabla u\|_{L^4}+\|\nabla H\|_{L^4}\big)^{\frac12}
\big(\|\nabla u_t\|_{L^2}+\|\nabla H_t\|_{L^2}\big) \notag\\
& \le C\big(\|\sqrt{\rho}\theta_t\|_{L^2}+\|\nabla\theta_t\|_{L^2}\big)
\big(\|\sqrt{\rho}\theta_t\|_{L^2}+\|\nabla\theta\|_{L^2}\big)^{\frac12}
\big(\|\nabla u_t\|_{L^2}+\|\nabla H_t\|_{L^2}\big) \notag \\
& \le \frac{\kappa}{12}\|\nabla \theta_t\|_{L^2}^2
+C\big(\|\sqrt{\rho}\theta_t\|_{L^2}^2
+\|\nabla\theta\|_{L^2}^2+1\big)
\big(\|\nabla u_t\|_{L^2}^2+\|\nabla H_t\|_{L^2}^2\big)
+\|\sqrt{\rho}\theta_t\|_{L^2}^2.
\end{align}
Therefore, inserting the above estimates on $L_1$--$L_3$ into \eqref{3.8.8} and combining \eqref{3.8.5} and \eqref{3.8.7}, we find that
\begin{align}\label{z3.8.11}
& \frac{d}{dt}\big(c_v\|\sqrt{\rho}\theta\|_{L^2}^2
+\kappa\|\nabla\theta\|_{L^2}^2+c_v\|\sqrt{\rho}\theta_t\|_{L^2}^2\big)
+\kappa\|\nabla\theta\|_{L^2}^2+c_v\|\sqrt{\rho}\theta_t\|_{L^2}^2
+\kappa\|\nabla\theta_t\|_{L^2}^2
\notag \\
&\le C\big(1+\|\nabla u_t\|_{L^2}^2+\|\nabla H_t\|_{L^2}^2\big)
\big(1+\|\sqrt{\rho}\theta\|_{L^2}^2+\|\nabla\theta\|_{L^2}^2
+\|\sqrt{\rho}\theta_t\|_{L^2}^2\big)+C\|\nabla^2\theta\|_{L^2}^2.
\end{align}

4. We deduce from \eqref{mhd}$_3$, the standard $L^2$-estimate of elliptic equations, \eqref{3.2}, \eqref{3.a2}, \eqref{2.4}, \eqref{3.5}, and \eqref{iAMSS12} that
\begin{align}\label{3.8.12}
 \|\nabla^2\theta\|_{L^2}^2
&\le C \big(\|\rho\theta_t\|_{L^2}^2+ \|\rho u\cdot\nabla\theta\|_{L^2}^2
+ \|\nabla u\|_{L^4}^4+ \|\nabla H\|_{L^4}^4\big) \notag \\
&\le C \big(\|\sqrt{\rho}\theta_t\|_{L^2}^2
+\|\rho u\|_{L^4}^2\|\nabla\theta\|_{L^4}^2
+\|\nabla u\|_{L^2}^2\|\nabla^2u\|_{L^2}^2
+\|\nabla H\|_{L^2}^2\|\nabla^2H\|_{L^2}^2\big) \notag \\
&\le C\|\sqrt{\rho}\theta_t\|_{L^2}^2
+C\|\nabla\theta\|_{L^2}\|\nabla^2\theta\|_{L^2}+C \notag \\
&\le C\|\sqrt{\rho}\theta_t\|_{L^2}^2
+\frac12\|\nabla^2\theta\|_{L^2}^2+C\|\nabla\theta\|_{L^2}^2+C,
\end{align}
which leads to
\begin{align}\label{3.8.13}
\|\nabla^2\theta\|_{L^2}^2
\leq C\|\sqrt{\rho}\theta_t\|_{L^2}^2+C\|\nabla\theta\|_{L^2}^2+C.
\end{align}
Hence, we derive from \eqref{z3.8.11} and \eqref{3.8.13} that
\begin{align}\label{zz3.8.11}
& \frac{d}{dt}\big(c_v\|\sqrt{\rho}\theta\|_{L^2}^2
+\kappa\|\nabla\theta\|_{L^2}^2+c_v\|\sqrt{\rho}\theta_t\|_{L^2}^2
\big)
+\kappa\|\nabla\theta\|_{L^2}^2+c_v\|\sqrt{\rho}\theta_t\|_{L^2}^2
+\kappa\|\nabla \theta_t\|_{L^2}^2
\notag \\
&\le C\big(1+\|\nabla u_t\|_{L^2}^2+\|\nabla H_t\|_{L^2}^2\big)
\big(1+\|\sqrt{\rho}\theta\|_{L^2}^2+\|\nabla\theta\|_{L^2}^2
+\|\sqrt{\rho}\theta_t\|_{L^2}^2\big).
\end{align}
Moreover, it follows from \eqref{mhd}$_3$, \eqref{c.c}, and \eqref{3.53} that
\begin{align*}
\int \rho\theta_t^2(x,0)dx
& \le\lim_{t \rightarrow 0}\sup\int\rho^{-1}\Big[\frac{\kappa}{c_v}\Delta\theta
+\frac{\mu}{2c_v}|\nabla u+(\nabla u)^{tr}|^{2}+\frac{\nu}{c_v}(\curl H)^{2}
-\rho u\cdot\nabla\theta\Big]^2dx \\
& \le C\|g_2\|_{L^2}^2+C\|\sqrt{\rho}u(0)\|_{L^\infty}^2\|\nabla\theta(0)\|_{L^2}^2 \le C,
\end{align*}
which combined with \eqref{zz3.8.11}, Gronwall's inequality, and \eqref{gj6} gives that
\begin{align}\label{3.8.11''}
&\sup\limits_{0\le t\le T}\big(\|\sqrt{\rho}\theta\|_{L^2}^2
+\|\nabla \theta\|_{L^2}^2+\|\sqrt{\rho}\theta_t\|_{L^2}^2\big)
+\int_0^T\big(\|\nabla\theta\|_{L^2}^2+\|\sqrt{\rho}\theta_t\|_{L^2}^2
+\|\nabla\theta_t\|_{L^2}^2\big)dt\le C.
\end{align}

5. One gets from \eqref{3.8.13} and \eqref{3.8.11''} that
\begin{align}\label{3.8.14}
\sup\limits_{0\le t\le T}
\|\nabla^2\theta\|_{L^2}^2+\int_0^T\|\nabla^2\theta\|_{L^2}^2dt\le C.
\end{align}
The standard $L^q$-estimate of elliptic equations together with \eqref{mhd}$_3$, \eqref{3.2}, H{\"o}lder's inequality, \eqref{local1}, \eqref{2.4}, \eqref{3.8.11''}, \eqref{3.4}, \eqref{3.8.14}, and \eqref{iAMSS12} yields that
\begin{align}\label{cz}
 \|\nabla^2\theta\|_{L^q}
& \le C \left(\|\rho \theta_t\|_{L^q}+\|\rho u\cdot\nabla\theta\|_{L^q} + \|\nabla u\|_{L^{2q}}^2+\|\nabla H\|_{L^{2q}}^2\right) \notag \\
& \le C\|\sqrt{\rho}\theta_t\|_{L^2}^{\frac{2(q-1)}{q^2-2}}
\|\sqrt{\rho}\theta_t\|_{L^{q^2}}^{\frac{q^2-2q}{q^2-2}}
+C\|\rho u\|_{L^{2q}}\|\nabla\theta\|_{L^{2q}} \notag \\
& \quad +C\|\nabla u\|_{L^2}^{\frac2q}\|\nabla^2u\|_{L^2}^{\frac{2q-2}{q}}
+C\|\nabla H\|_{L^2}^{\frac2q}\|\nabla^2H\|_{L^2}^{\frac{2q-2}{q}}
\notag \\
& \le C\|\sqrt{\rho}\theta_t\|_{L^2}^{\frac{2(q-1)}{q^2-2}}
\big(\|\sqrt{\rho}\theta_t\|_{L^2}+\|\nabla\theta_t\|_{L^2}\big)^{\frac{q^2-2q}{q^2-2}}
+C\big(\|\sqrt{\rho}u\|_{L^2}+\|\nabla u\|_{L^2}\big)
\|\nabla \theta\|_{L^2}^{\frac1q}\|\nabla^2\theta\|_{L^2}^{\frac{q-1}{q}}\notag \\
& \quad +C\|\nabla u\|_{L^2}^{\frac2q}\|\nabla^2u\|_{L^2}^{\frac{2q-2}{q}}
+C\|\nabla H\|_{L^2}^{\frac2q}\|\nabla^2H\|_{L^2}^{\frac{2q-2}{q}} \notag \\
& \le C\|\nabla\theta_t\|_{L^2}^{\frac{q^2-2q}{q^2-2}}+C,
\end{align}
which together with Young's inequality and \eqref{3.8.11''} indicates that
\begin{align} \label{3.8.16}
\int_0^T\Big(\|\nabla^2\theta\|_{L^q}^{\frac{q+1}{q}}
+\|\nabla^2\theta\|_{L^q}^2\Big)dt
\leq C\int_0^T\big(1+\|\nabla\theta_t\|_{L^2}^2\big)dt
\le C.
\end{align}
The proof of Lemma \ref{lemma 3.8} is finished.   \hfill $\Box$

\section{Proof of Theorem \ref{t1}}\label{sec4}

By Lemma \ref{lem21}, there exists a $T_{*}>0$ such that the problem \eqref{mhd}--\eqref{n4} has a unique local strong solution $(\rho,u,\theta,H)$ on $\mathbb{R}^2\times(0,T_{*}]$. We plan to extend the local solution to all time.

Set
\begin{equation}\label{20.1}
T^{*}=\sup \{T~|~(\rho,u,\theta,H)\ \text{is a strong solution on}\ \mathbb{R}^2\times(0,T]\}.
\end{equation}
First, for $T_*<T\leq T^{*}$ with $T$ finite, one deduces from \eqref{3.5}, \eqref{gj10'}, and \eqref{wdgj} that
\begin{equation}\label{20.2}
\nabla u, \nabla\theta, \nabla H\in C([0,T];H^1),
\end{equation}
where one has used the following fact (see \cite[Theorem 4, p. 304]{E2010})
\begin{equation}\label{20.0}
\|f\|_{C([0,T];H^1)}
\leq C(T) (\|f\|_{L^2(0,T;H^2)}+\|f_t\|_{L^2(0,T;L^2)}).
\end{equation}
Moreover, it follows from \eqref{igj1} and \eqref{gj8} that
\begin{equation}\label{20.3}
\rho\in C([0,T];L^1\cap H^1\cap W^{1,q}).
\end{equation}

Owing to \eqref{3.2} and \eqref{gj6}, we deduce that
\begin{align*}
\rho u_t=\sqrt{\rho}\cdot\sqrt{\rho}u_t\in L^2(0,T;L^2).
\end{align*}
Noting that
\begin{align*}
\rho_tu=-(u\cdot\nabla\rho)u
=-(u\bar{x}^{-\frac{a}{2}}\cdot\nabla\rho\bar{x}^a)u\bar{x}^{-\frac{a}{2}},
\end{align*}
which together with H{\"o}lder's inequality, \eqref{gj8}, and \eqref{3.a2} implies that
\begin{align*}
\rho_tu\in L^\infty(0,T;L^2).
\end{align*}
Thus, we arrive at
\begin{align}\label{20.4}
(\rho u)_t=\rho u_t+\rho_tu\in L^2(0,T;L^2).
\end{align}
From \eqref{3.3} and \eqref{3.4}, we have
\begin{align*}
\rho u=\sqrt{\rho}\cdot\sqrt{\rho}u\in L^\infty(0,T;L^2),
\end{align*}
which combined with \eqref{20.4} yields
\begin{align}\label{20.5}
\rho u\in C([0,T];L^2).
\end{align}
Similarly, one has
\begin{align}\label{20.6}
\rho\theta, H\in C([0,T];L^2).
\end{align}

Finally, if $T^{*}<\infty$, it follows from \eqref{20.2}, \eqref{20.3}, \eqref{20.5}, and \eqref{20.6} that
$$(\rho,u,\theta,H)(x,T^*)=\lim_{t\rightarrow T^*}(\rho, u,\theta,H)(x,t)$$
satisfies the initial condition \eqref{A} at $t=T^*$.
Furthermore, standard arguments yield that $(\rho\dot{u},\rho\dot{\theta})\in C([0,T^\ast];L^2)$,
which
implies
\begin{align*}
(\rho\dot{u},\rho\dot{\theta})(x,T^\ast)
=\lim_{t\rightarrow
T^\ast}(\rho\dot{u},\rho\dot{\theta})(x,t)\in L^2.
\end{align*}
Hence,
\begin{equation*}
\begin{split}
\begin{cases}
(-\mu\Delta u+\nabla P-H\cdot\nabla H)|_{t=T^\ast}
=\sqrt{\rho}(x,T^\ast)g_1(x)
,\\
(\kappa\Delta\theta
+\frac{\mu}{2}|\nabla u+(\nabla u)^{tr}|^2+\nu(\curl H)^2)|_{t=T^\ast}
=\sqrt{\rho}(x,T^\ast)g_2(x),
\end{cases}
\end{split}
\end{equation*}
with
\begin{align*}
g_1(x)\triangleq
\begin{cases}
\rho^{-\frac12}(x,T^\ast)(\rho\dot{u})(x,T^\ast),&
\mbox{for}~~x\in\{x|\rho(x,T^\ast)>0\},\\
0,&\mbox{for}~~x\in\{x|\rho(x,T^\ast)=0\},
\end{cases}
\end{align*}
and
\begin{align*}
g_2(x)\triangleq
\begin{cases}
\rho^{-\frac12}(x,T^\ast)(\rho\dot{\theta})(x,T^\ast),&
\mbox{for}~~x\in\{x|\rho(x,T^\ast)>0\},\\
0,&\mbox{for}~~x\in\{x|\rho(x,T^\ast)=0\},
\end{cases}
\end{align*}
satisfying $g_1,g_2\in L^2(\mathbb{R}^2)$ due to \eqref{gj6}, \eqref{ZZ}, \eqref{3.5}, and \eqref{3.8.11''}. So $(\rho, u,\theta,H)$ satisfies the compatibility condition \eqref{c.c} at $t=T^*$. Thus, taking $(\rho, u,\theta,H)(x,T^*)$ as the initial data, Lemma \ref{lem21} implies that one can extend the strong solutions beyond $T^*$. This contradicts the assumption of $T^*$ in \eqref{20.1}. Furthermore, the estimates as those in \eqref{1.10} follow from Lemmas \ref{l3.00}--\ref{lemma 3.8}. \eqref{l1.2} and \eqref{x5} follow from \eqref{p1} and \eqref{x4}, respectively. This completes the proof of Theorem \ref{t1}.  \hfill $\Box$

\end{document}